\documentclass[UTF-8,reqno]{amsart}
  
 \usepackage{enumerate}
\usepackage{amssymb,url,color, booktabs}
\usepackage{mathrsfs}
\usepackage{color}
\definecolor{MyDarkBlue}{cmyk}{0.8,0.3,0.8,0.4}
\definecolor{yellow}{rgb}{0.99,0.99,0.70}
\definecolor{white}{rgb}{1.0,1.0,1.0}
\definecolor{black}{rgb}{0.00,0.00,0.00}

\newcommand{\blue}{\color{black}}

\newcommand{\green}{\color{black}}

\usepackage[colorlinks=true]{hyperref}
\hypersetup{
    linkcolor=blue,          
    citecolor=red,        
    filecolor=blue,      
    urlcolor=cyan
}

\numberwithin{equation}{section}

\newcommand{\be}{\begin{eqnarray}}
\newcommand{\ee}{\end{eqnarray}}
\newcommand{\ce}{\begin{eqnarray*}}
\newcommand{\de}{\end{eqnarray*}}
\newtheorem{theorem}{Theorem}[section]
\newtheorem{lemma}[theorem]{Lemma}
\newtheorem{remark}[theorem]{Remark}
\newtheorem{definition}[theorem]{Definition}
\newtheorem{proposition}[theorem]{Proposition}
\newtheorem{Examples}[theorem]{Example}
\newtheorem{corollary}[theorem]{Corollary}

\def\eps{\varepsilon}

\def\e{\mathrm{e}}

\def\p{\partial}

\def\[{{\Big[}}
\def\]{{\Big]}}
\def\<{{\langle}}
\def\>{{\rangle}}
\def\({{\Big(}}
\def\){{\Big)}}

\def\bx{{\mathbf{x}}}

\def\dif{{\mathord{{\rm d}}}}

\def\no{\nonumber}
\def\={&\!\!=\!\!&}

\def\cB{{\mathcal B}}
\def\cC{{\mathcal C}}
\def\cD{{\mathcal D}}

\def\cL{{\mathcal L}}

\def\sS(\R^d){{\mathcal S}}

\def\mD{{\mathbb D}}
\def\mE{{\mathbb E}}

\def\mI{{\mathbb I}}

\def\mN{{\mathbb N}}

\def\mP{{\mathbb P}}

\def\mR{{\mathbb R}}

\def\1{{\mathbf{1}}}

\def\sA{{\mathscr A}}
\def\sB{{\mathscr B}}
\def\sC{{\mathscr C}}

\def\sK{{\mathscr K}}
\def\sL{{\mathscr L}}

\def\sS{{\mathscr S}}

\def\geq{\geqslant}
\def\leq{\leqslant}
\def\ge{\geqslant}
\def\le{\leqslant}

\def\div{\mathord{{\rm div}}}

\def\eps{\varepsilon}

\def\e{\mathrm{e}}

\def\p{\partial}

\def\[{{\Big[}}
\def\]{{\Big]}}
\def\<{{\langle}}
\def\>{{\rangle}}
\def\({{\Big(}}
\def\){{\Big)}}

\def\bx{{\mathbf{x}}}

\def\dif{{\mathord{{\rm d}}}}

\def\no{\nonumber}
\def\={&\!\!=\!\!&}
\def\bt{\begin{theorem}}
\def\et{\end{theorem}}
\def\bl{\begin{lemma}}
\def\el{\end{lemma}}
\def\br{\begin{remark}}
\def\er{\end{remark}}
\def\bx{\begin{Examples}}
\def\ex{\end{Examples}}
\def\bd{\begin{definition}}
\def\ed{\end{definition}}
\def\bp{\begin{proposition}}
\def\ep{\end{proposition}}
\def\bc{\begin{corollary}}
\def\ec{\end{corollary}}
\def\bpf{\begin{proof}}
\def\epf{\end{proof}}

\def\wt{\widetilde}

\def\geq{\geqslant}
\def\leq{\leqslant}
\def\ge{\geqslant}
\def\le{\leqslant}

\def\div{\mathord{{\rm div}}}

 \def\R{\mathbb R}
 \def\R{\mathbb R}    
\def\N{\mathbb N}  
   
\def\<{\langle} \def\>{\rangle}

 \def\beq{\begin{equation}}  
 
\def\e{\text{\rm{e}}}

\allowdisplaybreaks

\begin{document}

\title[Heat kernel of supercritical SDEs with unbounded drifts]{Heat kernel of supercritical nonlocal operators with unbounded drifts}

\author{St\'ephane Menozzi and Xicheng Zhang}

\address{Stephane Menozzi: Laboratoire de Mod\'elisation Math\'ematique d'Evry (LaMME), UMR CNRS 8071,
Universit\'e d'Evry Val d'Essonne (Universit\'e Paris Saclay), 23 Boulevard de France 91037 Evry, France and Laboratory of Stochastic Analysis, HSE,
Pokrovsky Blvd, 11, Moscow, Russian Federation\\ Email: stephane.menozzi@univ-evry.fr }

\address{Xicheng Zhang:
School of Mathematics and Statistics, Wuhan University,
Wuhan, Hubei 430072, P.R.China\\
Email: XichengZhang@gmail.com
 }

\thanks{
X. Zhang is supported by NNSF of China (Nos. 11731009, 12131019) and the DFG through the CRC 1283 
``Taming uncertainty and profiting from randomness and low regularity in analysis, stochastics and their applications''.\\
\indent S. Menozzi has been funded by the Russian Science Foundation project (project No. 20-11-20119).
}
\begin{abstract}
	
Let $\alpha\in(0,2)$ and $d\in\mN$. Consider the following SDE in $\mR^d$:
$$
\dif X_t=b(t,X_t)\dif t+a(t,X_{t-})\dif L^{(\alpha)}_t,\ \ X_0=x,
$$
where $L^{(\alpha)}$ is a $d$-dimensional rotationally invariant  $\alpha$-stable process, $b:\mR_+\times\mR^d\to\mR^d$ and $a:\mR_+\times\mR^d\to\mR^d\otimes\mR^d$ are  H\"older continuous functions in space, with respective order $\beta,\gamma\in (0,1)$ such that $(\beta\wedge \gamma)+\alpha>1$,  uniformly in $t$.  Here $b$ may be unbounded.
When $a$ is bounded and uniformly elliptic, we show that the unique solution $X_t(x)$ of the above SDE admits a continuous density, 
which enjoys sharp two-sided estimates. We also establish sharp upper-bound for the logarithmic derivative. 
{ In particular, we cover the whole \textit{supercritical} range $\alpha\in (0,1) $.}
Our proof is based on \textit{ad hoc} parametrix expansions and probabilistic techniques. 
	
\end{abstract}

\maketitle

\section{Introduction}

Throughout this paper we fix $\alpha\in(0,2)$. Let $L^{(\alpha)}$ be a $d$-dimensional rotationally invariant $\alpha$-stable process. 
We consider the following stochastic differential equation:
\begin{align}\label{SDE0}
\dif X_t=b(t,X_t)\dif t+a(t,X_{t-})\dif L^{(\alpha)}_t,
\end{align}
where $b:\mR_+\times\mR^d\to\mR^d$ and $a:\mR_+\times\mR^d\to\mR^d\otimes\mR^d$ are Borel measurable functions and satisfy
that for some $\beta\in((1-\alpha)^+,1]$ and $\kappa_0\geq 1$,
\begin{align}\label{BB0}
|b(t,0)|\le \kappa_0,\ |b(t,x)-b(t,y)|\leq\kappa_0(|x-y|^\beta\vee|x-y|),\tag{{\bf H$^\beta_b$}}
\end{align}
and 
for some $\gamma\in((1-\alpha)^+,1]$ and $\kappa_1\geq 1$,
\begin{align}\label{AA0}
\kappa_1^{-1}\mI\leq (aa^*)(t,x)\leq\kappa_1\mI,\ \ |a(t,x)-a(t,y)|\leq\kappa_1|x-y|^{\gamma}, \tag{{\bf H$^{\gamma}_a$}}
\end{align}
where $a^*$ stands for the transpose of $a$ and $\mI$ is the identity matrix.
Under {\bf (H$^\beta_b$)} and {\bf (H$^{\gamma}_a$)}, it is well known that 
for each $(s,x)\in\mR_+\times\mR^d$,
there is a unique weak solution 
$X_{s,t}(x)$ to SDE \eqref{SDE0} starting from $x$ at time $s$ (see e.g. \cite[Theorem 1.1]{CZZ}),
and the generator of SDE \eqref{SDE0} writes as
\begin{align}
\sL_{s} f(x):=\tfrac12\cL_{s} f(x)+b(s,x)\cdot\nabla f(x),
\end{align}
where $\cL_{s}$ is given by
\begin{align}\label{GEN_STABLE_LIKE}
\cL_{s} f(x)=\int_{\mR^d}\delta^{(2)}_f(x;a(s,x)z)\frac{\dif z}{|z|^{d+\alpha}}=\int_{\mR^d}\delta^{(2)}_f(x;z)\frac{\kappa(s,x,z)}{|z|^{d+\alpha}}\dif z
\end{align}
with 
\begin{equation}\label{DEF_DELTA_2}
\delta^{(2)}_f(x;z):=f(x+z)+f(x-z)-2f(x)
\end{equation}
and
\begin{align}\label{KA1}
\kappa(s,x,z):=\det(a^{-1}(s,x))(|z|/|a^{-1}(s,x)z|)^{d+\alpha}.
\end{align}
Clearly, by \eqref{AA0} we have for some $\bar\kappa_1\geq 1$,
\begin{align}\label{AA2}
\bar\kappa^{-1}_1\leq \kappa(s,x,z)\leq\bar\kappa_1,\ \ |\kappa(s,x,z)-\kappa(s,y,z)|\leq\bar\kappa_1|x-y|^{\gamma}.
\end{align}
The operator $\sL_{s}$ is called \textit{supercritical}   for $\alpha\in(0,1)$ since in this case, the drift term plays a dominant role. Namely, from the self-similarity properties of the driving process $L^{(\alpha)}$  in \eqref{SDE0}, it holds that for any $s>0$, $L_s^{(\alpha)} \overset{({\rm law})}{=} s^{1/\alpha}L_1^{(\alpha)} $ and for  $s\in (0,1),\ \alpha\in (0,1)$, $s^{1/\alpha}< s $. This precisely means that the fluctuations induced by the noise are smaller than the typical order of the drift term in \eqref{SDE0}. 
For $\alpha\in(1,2)$, the converse phenomenon happens. Since for $s\in (0,1) $, $s^{1/\alpha}>s $, the fluctuations of the noise prevail in the SDE. From the operator viewpoint, $\cL_{s}$ plays a dominant role and we say that $\sL_{s}$ is \textit{subcritical}. For the remaining case $\alpha=1$, the noise and drift both have the same typical order and the operator $\sL_{s}$ is called \textit{critical}.
Note that for $\alpha\in(0,1)$, since $z\mapsto\kappa(s,x,z)$ is symmetric, we have
$$
\cL_{s} f(x)=2\int_{\mR^d}\delta^{(1)}_f(x;z)\frac{\kappa(s,x,z)\dif z}{|z|^{d+\alpha}},
$$
where
$$
\delta^{(1)}_f(x;z):=f(x+z)-f(x).
$$

\smallskip

Let us now indicate that there is a quite large literature concerning stable driven SDEs. We can first mention the seminal work of Kolokoltsov \cite{kolo:00} \textcolor{black}{from which one can derive} that  for an SDE driven by a symmetric stable process with smooth non-degenerate spectral measure, Lipschitz non-degenerate diffusion coefficient and non trivial Lipschitz bounded drifts when $\alpha>1 $,  two sided estimates for the density of the type:
\begin{equation}\label{FIRST_TW}
p(s,x,t,y)\asymp_C (t-s)^{-d/\alpha} \left(1+\frac{|x-y|}{(t-s)^{1/\alpha}}\right)^{-(d+\alpha)},
\end{equation}
where $C\ge 1$ depends on the non-degeneracy and Lipschitz constants of the coefficients and the final considered time horizon $T$.
Here and below, $Q_1\asymp_C Q_2$ means that $C^{-1}Q_2\leq Q_1\leq C Q_2$.

\smallskip

Going to weaker regularity of the coefficients in \eqref{SDE0} then first leads to investigate the well-posedness of the martingale problem associated with the formal generator associated with the dynamics \eqref{SDE0}. In \cite{Ba-Ch},
Bass and Chen showed the weak well-posedness for SDE \eqref{SDE0} when $a$ is only continuous and uniformly elliptic, 
$b$ is Lipschitz and $L^{(\alpha)}_t$ is cylindrical $\alpha$-stable process.
In the subcritical case we can mention the work by Mikulevicius and Pragarauskas \cite{miku:prag:14} 
who derived that weak uniqueness holds for equation \eqref{SDE0} for bounded H\"older coefficients when  $ \alpha\ge 1$ 
and a non degenerate $a$. The martingale problem was in their framework studied from some related Schauder estimates established on the associated Integro Partial Differential Equation (IPDE). We can also refer to \cite{chau:meno:prio:20} for \textcolor{black}{parabolic} Schauder estimates in  the super-critical \textcolor{black}{stable} case and \textcolor{black}{to the work by K\"uhn \cite{kuhn:20} for interior elliptic Schauder estimates for a larger class of L\'evy operators satisfying the Hartman-Wintner growth condition and suitable $L^1$ gradient estimates for the corresponding heat kernel}.

\smallskip

In the super-critical case, the well-posedness of the martingale problem was recently investigated by 
Kulik \textit{et al.} \cite{knop:kuli:18}, \cite{kuli:19} (see also \cite{CZZ}). 
In \cite{kuli:19}, 
the authors consider SDEs of type \eqref{SDE0} with bounded H\"older drift and non-degenerate scalar diffusion coefficients under the natural condition $\alpha+\beta>1 $\footnote{it is indeed well known that from the seminal work of Tanaka \textit{et al.} \cite{tana:tsuc:wata:74} that weak uniqueness may fail if this condition is not met.} and obtain the existence of the heat kernel, 
a corresponding two-sided estimate of the form \eqref{FIRST_TW} as well as some estimates corresponding to the time derivative through parametrix type expansions. Let us emphasize that in the super-critical regime, the time derivative of the heat kernel roughly typically behaves as $t^{-1}$ at time $t$ whereas the spatial gradient is then more singular, as it is expected to have typical behavior of order $t^{-1/\alpha}\ge t^{-1} $ for $t\in (0,1] $.

\smallskip

\textcolor{black}{Concerning other results related to  stable heat kernel estimates we can refer e.g. 
for driftless operators
of the form \eqref{GEN_STABLE_LIKE} to \cite{Ch-Zh}
  \cite{Ch-Zh1} in which the authors consider non-symmetric functions $\kappa$ with possible unbounded dependence in the jump variable in the second paper. Therein, two-sided heat kernel and gradient estimates are derived, and for a bounded, even $\kappa:=\kappa(x) $ the gradient estimate has been shown to be sharp in \cite{DZ}. We can also quote \cite{CHXZ} for similar results when the operator also has a diffusion part, or \cite{Ch-Zh20} for more general absolutely continuous L\'evy measures satisfying suitable scaling properties.
}
  
\textcolor{black}
{We insist that, beyond the symmetry and the absolute continuity condition of the L\'evy measure w.r.t. the one of the isotropic stable process, the behavior of the stable densities or stable driven SDEs can be very different from the above description. This can be 	already seen from the seminal work of Watanabe \cite{wata:07} for stable processes and from the recent work \cite{knop:kuli:schi:20} for SDEs, see Section 4 therein, and in particular Example 4.2 which relates to a cylindrically stable driven SDE whose density turns out to be unbounded.
}

\textcolor{black}{
For non zero drifts and rather general, possibly non-symmetric, stable L\'evy measures, estimates of the heat kernel in Besov norms were obtained in \cite{De-Fo}.
We also refer to the work \cite{kuli:pesz:prio:20}  and \cite{CHZ} for gradient estimates on the semigroup associated 
with additive and multiplicative cylindrical noises in \eqref{SDE0} respectively, for $\alpha\in (0,2) $ and appropriate assumptions on the drift. Still in the cylindrical case, we can quote the papers by Kulczycki, Ryznar \textit{et al.} \cite{kulc:ryzn:szto:21} for semigroup estimates in the stable super-critical case, \cite{kulc:ryzn:19}, \cite{kulc:kuli:ryzn:21} for semigroup and diagonal estimates on the density for more general driving processes whose characteristic exponents satisfy upper and lower scaling conditions.
}

\smallskip

In the current work, \textcolor{black}{we stick to operators of the form \eqref{GEN_STABLE_LIKE}, but anyhow face two difficulties:} we want to establish \textcolor{black}{density and} gradient estimates for all $\alpha\in (0,2) $ and for unbounded drifts. It is known, and somehow intuitive, that for unbounded drifts the heat kernel bounds must reflect somehow the transport induced by the drift. This was for instance observed for a Lipschitz drift in \cite{dela:meno:10} for degenerate Kolmogorov SDEs which can be viewed as ODEs perturbed on some components by a Brownian noise propagating through the whole chain thanks to a weak type H\"ormander condition on the drift. 
Before going further let us also mention the work by Huang \cite{huan:15} 
which establishes two-sided estimates for stable driven SDE with unbounded Lipschitz drift and $\alpha\in(1,2)$, 
which then read as: for $s<t$ and $x,y\in\mR^d$,
\begin{equation}
\label{SECOND_TW}
p(s,x,t,y)\asymp_C (t-s)^{-d/\alpha} \left(1+\frac{|\theta_{s,t}(x)-y|}{(t-s)^{1/\alpha}}\right)^{-(d+\alpha)},
\end{equation}
where $\theta_{s,t}(x) $ denotes the flow associated to the drift in \eqref{SDE0}. Namely, 
$$
\dot \theta_{s,t}(x)=b(t,\theta_{s,t}(x)), \theta_{s,s}(x)=x.
$$

\smallskip

In the non-degenerate Brownian case, the type of heat kernel estimate { in} \cite{dela:meno:10} has recently been extended to drifts satisfying a linear-growth without \textit{a priori} smoothness assumptions on the drift, see \cite{meno:pesc:zhan:20}. In the quoted work, under additional H\"older continuity of the drift  for the second derivatives, the estimates also extend to the derivatives up to order two with the corresponding additional parabolic singularity.

\textcolor{black}{All the previously quoted works use the \textit{parametrix} method. Initiated by E.E. Levi \cite{levi:07},  it has become a classical tool for the study of heat kernels of operators with non constant coefficients, see e.g. \cite{Friedman}, \cite{mcke:sing:67} in the diffusive case or \textcolor{black}{\cite{koch:89} in a much larger framework}.  The idea of the method is that, in \textit{small} time, assuming that the heat kernel of the operator with variable coefficients exists, it should \textit{somehow} behave as a well understood one corresponding to an operator with \textit{constant} coefficient. The difference between the two densities is then controlled through a Duhamel type formula.
In its usual formulation this approach is designed for operators with spatially bounded coefficients. Assume for a while that, additionally to the stated assumptions,  the coefficients in \eqref{SDE0} are smooth and bounded. It is then known that the density of $(X_t)_{t\ge s}$ exists for $t>s$. Denote it by $p(s,x,t,\cdot)$ for the starting point $x$ at time $s$. The \textit{parametrix} is a zero order approximation of the density $p(s,x,t,\cdot) $ to be estimated and for which one precisely knows the behavior, e.g. Gaussian in \cite{Friedman}, stable isotropic with possible shift in \cite{kolo:00}. The most common choice consists in considering, for a fixed terminal point $y$,  as parametrix the density of the process with dynamics
\begin{equation}\label{PROXY_1}
\dif \tilde X_t^y=b(t,y)\dif t+a(t,y)\dif L^{(\alpha)}_t,\ t\ge s,\ \tilde X_s^y=x.
\end{equation}
From the uniform ellipticity of $a$ and the specific form \eqref{GEN_STABLE_LIKE} it is possible to derive that 
$ \tilde X_t^y$ admits a density for $t>s$ and that for any $k\in \N,\ w\in \R^d$,
\begin{align}\label{GRAD_BD_INTRO}
|\nabla_x^k \tilde p^y(s,x,t,w)|\le C (t-s)^{-\frac{ k+d}\alpha }\Big(1+\frac{|w-(x+\int_s^t b(r,y)\dif r )|}{(t-s)^{\frac 1\alpha}}\Big)^{-(d+\alpha+k)},
\end{align}
see e.g. \cite{Be}, \cite{Bo-Ja} and Section \ref{SEC_DENS_EST} below. The difference between the density of the SDE and the \textit{proxy} one is then investigated through the Duhamel formula which formally follows from the Kolmogorov equations satisfied by the initial density and the \textit{parametrix} in \eqref{PROXY_1}. Namely, $\tilde p^y(\cdot,\cdot,t,w) $ is a classical solution of:
\begin{equation}\label{KOLMO_FREEZE_INTRO}
\p_s \tilde p^{y}(s,x,t,w)+{\tilde \sL}^{y}_s{\tilde p}^{y}(s,\cdot,t,w)(x)=0,\ (s,x)\in [0,t)\times \R^d,
\end{equation}
where
$$
\tilde \sL^{y}_{s} f(x):=\tfrac 12\tilde \cL^{y}_{s} f(x)+b(s, y)\cdot\nabla f(x),\
\tilde \cL^{y}_{s} f(x)=\int_{\mR^d}\delta^{(2)}_f(x;z)\frac{\kappa(s,y,z)}{|z|^{d+\alpha}}\dif z
$$
with
$$
\kappa(s,y,z):=\frac{\det(a^{-1}(s,y)|z|^{d+\alpha}}{|a^{-1}(s,y)z|^{d+\alpha}}.
$$
And $\tilde p^{y}(s,\cdot,t,w) \underset{s\uparrow t}{\longrightarrow}  \delta_w(\cdot)$, weakly identifying the density with the induced measure. Equation \eqref{KOLMO_FREEZE_INTRO} is the \textit{backward} Kolmogorov equation. If the density of the SDE is now a classical solution of the forward Kolmogorov equation
\begin{align}\label{KOLMO_INTRO}
\p_t  p(s,x,t,w)&-({ \sL}_t)^*p(s,x,t,\cdot)(w)=0,\ (t,w)\in (s,+\infty]\times \R^d,\\
p(s,x,t,w) &\underset{t\downarrow s}{\longrightarrow}  \delta_x(\cdot),\notag
\end{align}
where $({ \sL}_t)^* $ denotes the adjoint of ${ \sL}_t $,
it formally follows that
\begin{align}
(p-\tilde p^y)(s,x,t,y)=&\int_s^t \dif r \partial_r\Big(\int_{\R^d}p(s,x,r,w) \tilde p^y(r,w,t,y) \dif w\Big)\notag\\
=&\int_s^t \dif r  \int_{\R^d} \partial_r p(s,x,r,w) \tilde p^y(r,w,t,y) \dif w\notag\\
&+\int_s^t \dif r \int_{\R^d}  p(s,x,r,w) \partial_r\tilde p^y(r,w,t,y) \dif w\notag\\
\underset{\eqref{KOLMO_FREEZE_INTRO}, \eqref{KOLMO_INTRO}}{=}
&\int_s^t \dif r \int_{ \R^d} p(s,x,r,w)({ \sL}_r-\tilde { \sL}_r^y)\tilde p^y (r,\cdot,t,y)(w)\dif w,
\label{PARAM_INTRO}
\end{align}
 provided all the previous computations can be justified, which can be delicate for the density of the SDE. The idea is then to repeat the approximation procedure in the integral of \eqref{PARAM_INTRO}, i.e. approximate $p(s,x,r,w)$ at order 0 with $\tilde p^w(s,x,r,w)$ and control again the difference through the associate integral. Iterating infinitely many times leads to the so-called parametrix series. To obtain some quantitative bounds from this procedure it is clear that one has to understand the behavior of the so-called \textit{parametrix kernel}:
 \begin{align}
q_0(r,w,t,y):=&(\sL_{r}-\tilde \sL^{y}_r)\tilde p^y (r,\cdot,t,y)(w)\notag\\
=&[\tfrac 12(\tilde \cL_r-\tilde \cL^{y}_{r}) +(b(r,w)-b(r, y))\cdot\nabla] \tilde p^y(r,\cdot,t,y)(w),\label{DEF_A_INTRO}
\end{align}
 and its smoothing properties. Let us here focus on the drift term. For a bounded Lipschitz drift and $\alpha\in [1,2)$, it is seen from \eqref{GRAD_BD_INTRO} that
 \begin{align}\label{DRIFT_KERNEL_BD_INTRO}
&|b(r,w)-b(r,y)||\nabla_w \tilde p^y(r,w,t,y)|\notag\\
\le& C |w-y|(t-r)^{-\frac{ 1+d}\alpha }\Big(1+\frac{|y-(w+\int_r^t b(u,y)\dif u)|}{(t-r)^{\frac 1\alpha}}\Big)^{-(d+\alpha+1)}\notag\\
\le& \tilde C (t-r)^{-\frac{d}\alpha }\Big(1+\frac{|y-w|}{(t-r)^{\frac 1\alpha}}\Big)^{-(d+\alpha)},
\end{align}
where for the last inequality we also used that, in small time and for a bounded drift $b $ and $\alpha\ge 1 $, the frozen drift $\int_r^t b(u,y)\dif u$ is negligible w.r.t. to the characteristic time scale, i.e. $|\int_r^t b(u,y)\dif u|\le \|b\|_\infty(t-r)\le \|b\|_\infty(t-r)^{\frac 1\alpha}  $ so that, up to a modification of $C$,
$$\Big(1+\frac{|y-(w+\int_r^t b(u,y)\dif u)|}{(t-r)^{\frac 1\alpha}}\Big)^{-(d+\alpha+1)}\le  C \Big(1+\frac{|y-w|}{(t-r)^{\frac 1\alpha}}\Big)^{-(d+\alpha+1)}.$$
We have thus bounded in \eqref{DRIFT_KERNEL_BD_INTRO} the drift part of the parametrix kernel by a quantity which serves as upper and lower bound for the isotropic stable density, see again \eqref{FIRST_TW}. Hence, the regularity of the drift allowed to absorb the time singularity induced by the spatial derivative of the \textit{parametrix}. Observe that if $\alpha\in (1,2) $ and since the drift is bounded, we could even have taken a  driftless parametrix. In that case, the contribution associated with the difference of the drifts, between the generator of the SDE and the one of the frozen process would write:
 \begin{align}\label{DRIFT_KERNEL_BD_INTRO_SUB_CRIT}
&|b(r,w)||\nabla_w \tilde p^y(r,w,t,y)|\notag\\
\le& C \|b\|_\infty(t-r)^{-\frac{ 1+d}\alpha }\Big(1+\frac{|y-w|}{(t-r)^{\frac 1\alpha}}\Big)^{-(d+\alpha+1)}\notag\\
\le& \tilde C (t-r)^{-\frac{d+1}\alpha }\Big(1+\frac{|y-w|}{(t-r)^{\frac 1\alpha}}\Big)^{-(d+\alpha)},
\end{align}
which is homogeneous to an isotropic stable density multiplied by a time integrable singularity. 
Having in mind that in \eqref{PARAM_INTRO} the parametrix kernel is integrated in time, 
this is sufficient to provide a smoothing effect in small time, i.e. the integral term in \eqref{PARAM_INTRO} will be negligible w.r.t. to the first order approximation, provided we manage to control the series expansion deriving from the iteration of the procedure. This requires to control  convolutions of the density of the parametrix and iterated convolutions of the parametrix kernel.}

\textcolor{black}{Turning now to the case $\alpha\in (0,1) $, even if $b$ is bounded, the last inequality of \eqref{DRIFT_KERNEL_BD_INTRO} fails precisely because $t^{\frac 1\alpha } \le t$, $t\le 1 $, i.e. the drift prevails w.r.t. characteristic time of the noise in small time. Thus, freezing at the terminal point does not provide a \textit{good} parametrix.  On the other hand, taking a driftless proxy is not an option either, since the associated bound for the drift part of the parametrix kernel would yield a non-integrable singularity in \eqref{DRIFT_KERNEL_BD_INTRO_SUB_CRIT}. We thus need to consider a proxy which can absorb the time singularity of the gradient through the difference of the drift coefficients in the operators. If $b$ is Lipschitz, this can e.g. be done through the backward flow associated with the differential system. Namely,
introducing for fixed final point $y\in \R^d$ and final time $t$,
\begin{equation}\label{DIFF_DYN}
\dot \theta_{t,r}(y)=b(r,\theta_{t,r}(y)), r\in [0,t],\ \theta_{t,t}(y)=y,
\end{equation}
and considering as parametrix the density of the process with dynamics
\begin{equation}\label{PROXY_2}
\dif \tilde X_r^{(t,y)}=b(r,\theta_{t,r}(y))\dif t+a(r,\theta_{t,r}(y))\dif L^{(\alpha)}_t,\ r\ge s,\ \tilde X_s^{(t,y)}=x,
\end{equation} 
one derives similarly to \eqref{GRAD_BD_INTRO} 
\begin{align}\label{GRAD_BD_INTRO_BIS}
|\nabla_w^k \tilde p^{(t,y)}(r,w,t,y)|\le &C (t-r)^{-\frac{ k+d}\alpha }\Big(1+\frac{|y-(w+\int_r^t b(u,\theta_{t,u}(y))\dif u)|}{(t-r)^{\frac 1\alpha}}\Big)^{-(d+\alpha+k)}\notag\\
\le & C (t-r)^{-\frac{ k+d}\alpha }\Big(1+\frac{|\theta_{t,r}(y)-w)|}{(t-r)^{\frac 1\alpha}}\Big)^{-(d+\alpha+k)},
\end{align}
exploiting the differential dynamics \eqref{DIFF_DYN} in backward time for the last inequality. 
Writing now the difference of the drift parts of the generators applied to the frozen density $\tilde p^{(t,y)}(\cdot,\cdot,t,y)$ yields:
 \begin{align*}
&|b(r,w)-b(r,\theta_{t,r}(y))||\nabla_w \tilde p^y(r,w,t,y)|\notag\\
\le& C |w-\theta_{t,r}(y)|(t-r)^{-\frac{ 1+d}\alpha }\Big(1+\frac{|\theta_{t,r}(y)-w|}{(t-r)^{\frac 1\alpha}}\Big)^{-(d+\alpha+1)}\notag\\
\le& \tilde C (t-r)^{-\frac{d}\alpha }\Big(1+\frac{|\theta_{t,r}(y)-w|}{(t-r)^{\frac 1\alpha}}\Big)^{-(d+\alpha)},
\end{align*}
and this choice then allows to absorb the time singularity similarly to the previously considered sub-critical case $\alpha\in [1,2) $, see \eqref{DRIFT_KERNEL_BD_INTRO}. For parametrix expansions and density estimates, this idea of using the backward flow associated with the first order term was already used in the diffusive setting in \cite{kona:meno:molc:10}, \cite{dela:meno:10} and for stable driven SDE in \cite{huan:15} for Lipschitz drifts with linear growth in the subcritical case and \cite{knop:kuli:18} for bounded Lipschitz drifts in the supercritical case. In this last setting, an underlying flow was also used in \cite{kuli:19} where the author addresses the case of a bounded, time-homogeneous H\"older drift coefficient\footnote{Let us mention that the boundedness of the drift coefficient in \cite{kuli:19} follows from the specific chosen zero order approximation. Namely the diffusion coefficient therein is not frozen along the flow but at the final spatial point, see e.g. the Proof of Lemma 4.1 in that reference.}. In this framework, there is no uniqueness to \eqref{DIFF_DYN} but only existence (thanks to the Peano theorem). Anyhow, for a $\beta$ H\"older drift, bounded or not, and a given associated Peano flow $\theta_{r,t}(y) $ solving \eqref{DIFF_DYN},
 \begin{align}\label{DRIFT_KERNEL_BD_INTRO_BD_FLOW}
&|b(r,w)-b(r,\theta_{t,r}(y))||\nabla_w \tilde p^y(r,w,t,y)|\notag\\
\le& C |w-\theta_{t,r}(y)|^\beta(t-r)^{-\frac{ 1+d}\alpha }\Big(1+\frac{|\theta_{t,r}(y)-w|}{(t-r)^{\frac 1\alpha}}\Big)^{-(d+\alpha+1)}\notag\\
\le& \tilde C (t-r)^{-(\frac{d}\alpha+\frac{1-\beta}{\alpha} )}\Big(1+\frac{|\theta_{t,r}(y)-w|}{(t-r)^{\frac 1\alpha}}\Big)^{-(d+\alpha)},
\end{align}
which yields an integrable singularity provided $\alpha+\beta>1 $ (which is the natural condition we assumed and already discussed above). This type of estimate is crucial for the parametrix approach to work in our framework.}

\textcolor{black}{In the following, in order to benefit from Lipschitz properties of an associate flow, which is e.g. needed to make convolutions of terms which have the form of the r.h.s. of \eqref{DRIFT_KERNEL_BD_INTRO_BD_FLOW} and which will naturally appear from the parametrix expansion, we will not use exactly \eqref{PROXY_2} as zero order approximation but replace therein $b $ by a suitable mollification observed along the associated backward (Lipschitz) flow. Rather naturally, the mollification parameter is chosen to correspond to the characteristic scaling time, see equation \eqref{ODE1} below}

\textcolor{black}{With respect to the various previously described steps, let us indicate that some  properties of mollified flows, that appear in our main results, will be stated at the end of the section. The bounds on the frozen densities and the related convolution estimates, which will appear when iterating the first order expansion in \eqref{PARAM_INTRO} will be discussed in Section \ref{SEC_PRELIM}.}\\

We will here somehow follow the main line of \cite{meno:pesc:zhan:20} but are faced with many additional difficulties. In particular, a common feature to both the Gaussian SDEs considered in \cite{meno:pesc:zhan:20} and the stable driven here is that we first need to establish the density and gradient estimates for smooth coefficients {\blue and a bounded drift, in order to justify that the previous expansions can indeed be performed (see in particular Lemma \ref{Le310}).} 
In the current strictly stable framework we cannot rely on the Malliavin calculus arguments of \cite{meno:pesc:zhan:20} because of integrability issues.
We thus establish here some direct bounds on the associated semi-group and its derivatives when the coefficients are smooth and \textcolor{black}{the drift bounded} which serve as a starting point to derive the estimates of Theorem \ref{Main}. 
This part is crucial and quite intricate (see Theorems \textcolor{black}{\ref{T8} and}  \ref{Th41} below).
Since the two-sided and gradient estimates hold independently of the smoothness of the coefficients {\blue and the boundedness of the drift}, 
we eventually conclude through compactness {\blue arguments (see Section \ref{FINAL_SEC})}.

\smallskip

To state our main result, we introduce the regularized flow associated with drift $b$.
For $\eps>0$, let $b_\eps(t,x):=b(t,\cdot)*\rho_\eps(x)$, where $\rho_\eps(x)=\eps^{-d}\rho(x/\eps)$
and $\rho$ is a smooth density function with support in the unit ball $B(0,1)$.
Note that under \eqref{BB0},
{\blue
\begin{align}\label{LL19}
|b_\eps(t,x)-b(t,x)|=&
\left|\int_{\R^d} \big(b(t,x-y)-b(t,x)\big)\rho_\eps(y)\dif y\right|\notag\\
\le &
\kappa_0   \int_{B(0,\eps)}( |y|^\beta+|y|)   \rho_\eps(y)\dif y\leq \kappa_0(\eps^\beta+\eps),
\end{align}
\textcolor{black}{since $\int_{B(0,\eps)}\rho_\eps(y)\dif y=1 $. Recalling now that $\nabla_x\underbrace{\int_{\R^d} \rho_\eps(x-y) \dif y}_{=1}=0 $},
\begin{align}\label{LL1}
|\nabla b_\eps(t,x)|
=&\left|\int_{\mR^d}(b(t,y)-b(t,x))\nabla\rho_\eps(x-y)\dif y\right|\notag\\
\leq&\kappa_0\int_{B(x,\eps)} \!\!(|x-y|^\beta+|x-y|)  \eps^{-(d+1)}|\nabla \rho(z)|_{z=\frac{x-y}\eps}\dif y
\notag\\
\le& \kappa_0(\eps^{\beta-1}+1)\|\nabla\rho\|_{L^1}.
\end{align}
}
In particular, since $\alpha+\beta>1$, for any $T>0$, there is a $C>0$ such that for any $0\leq s<t\leq T$,
\begin{align}\label{LL01}
\begin{split}
\int^t_s\|\nabla b_{|r-s|^{1/\alpha}}(r,\cdot)\|_\infty\dif r
&\leq C
\int^t_s\Big(|r-s|^{\frac{\beta-1}{\alpha}}+1\Big)\dif r
\leq C(t-s)^{\frac{\alpha+\beta-1}{\alpha}}.
\end{split}
\end{align}
Thus, for fixed $s>0$, the following ODE admits a unique solution $\theta_{s,t}(x)$: 
\begin{align}\label{ODE1}
\dot\theta_{s,t}=b_{|t-s|^{1/\alpha}}(t,\theta_{s,t}),\ \ \theta_{s,s}=x,\ \ \ t\geq 0.
\end{align}
Note that for $t>s$, $\theta_{s,t}(x)$ denotes the forward solution of the above ODE, while for $t<s$, 
it denotes the backward solution. We carefully mention that our main results will be stated w.r.t. to the flow $\theta $ in \eqref{ODE1} which is precisely associated with a mollified drift with parameter corresponding to the typical scale of the driving process of the SDE \eqref{SDE0} at the current considered time. 

For notational simplicity, we  introduce the following parameter set
\begin{equation}\label{DEF_THETA}
\Theta:=(\kappa_0,\kappa_1,d,\alpha,\beta,\gamma).
\end{equation}
We also denote for $T\in(0,\infty]$,
$$
\mD_T:=\{(s,x,t,y): 0\leq s<t<T, x,y\in\mR^d\}.
$$
We will frequently use from now on the notation $\lesssim$. For  two quantities $Q_1$ and $Q_2$, we mean by 
$Q_1\lesssim Q_2 $ that there exists $C:=C(T,\Theta)$  such that $Q_1\le CQ_2 $. Other possible dependencies for the constants will be explicitly specified.
Moreover, we also use the following notation
\begin{align}\label{DA0}
|\cD^{(\alpha)} f|(x):=\int_{\mR^d}\frac{|\delta^{(2)}_f(x;z)|}{|z|^{d+\alpha}}\dif z.
\end{align}

The aim of this paper is to show the following result.
\bt\label{Main}
Under {\bf (H$^\beta_b$)} and {\bf (H$^{\gamma}_a$)}, for each $0\leq s<t<\infty$ and $x\in\mR^d$, $X_{s,t}(x)$ admits a density $p(s,x,t,y)$ 
(called heat kernel of $\sL_{s}$) that is continuous as a function of $x,y$, and such that for each $t>0$ and $x,y\in\mR^d$ and Lebesgue almost all $s\in[0,t)$,
$$
\p_s p(s,x,t,y)=\sL_{s}p(s,\cdot,t,y)(x),\ \  p(s,x,t,\cdot)\rightarrow \delta_{\{x\}}(\cdot)\ {\rm \it weakly\ as}\ s\uparrow t,
$$
where $\delta_{\{x\}}(\dif y)$ denotes the Dirac measure concentrated at $x$.
Moreover, we have
\begin{enumerate}[(i)]
\item {\bf (Two-sides estimate)} For any $T>0$, there is a constant $C_1=C_1(T,\Theta)\ge 1$ such that for all $(s,x,t,y)\in\mD_T$,
\begin{align}\label{TW}
p(s,x,t,y)\asymp_{C_1} (t-s)((t-s)^{1/\alpha}+|\theta_{s,t}(x)-y|)^{-d-\alpha},
\end{align}
{ where $\theta_{s,t}(x)$ is defined by ODE \eqref{ODE1}}.
\item {\bf (Fractional derivative estimate)} For any $T>0$, there is a constant $C_2=C_2(T,\Theta)>0$ such that for all $(s,x,t,y)\in\mD_T$,
\begin{align}\label{FR}
|\cD^{(\alpha)}p(s,\cdot,t,y)|(x) \le {C_2} ((t-s)^{1/\alpha}+|\theta_{s,t}(x)-y|)^{-d-\alpha}.
\end{align}
\item {\bf (Gradient estimate in $x$)} For any $T>0$, there is a constant $C_3=C_3(T,\Theta)>0$ such that for all $(s,x,t,y)\in\mD_T$,
\begin{align}\label{GR0}
|\nabla_x\log p(s,x,t,y)|\le {C_3} (t-s)^{-1/\alpha}.
\end{align}
\end{enumerate}
\et

\br
\rm 
If $|b(t,x)-b(t,y)|\leq\kappa_0|x-y|^\beta$ for any $x,y\in\mR^d$ with $|x-y|\leq 1$ and $|b(0,t)|\le \kappa_0 $ for all $t\ge 0 $, then
{\bf (H$^\beta_b$)} holds. { In particular, for $c(x)$ being a bounded $\beta$-H\"older continuous function, $b(x):=x+c(x)$ satisfies {\bf (H$^\beta_b$)}.}
\er

\br\rm
For $\alpha\in[1,2)$, we can replace $\theta_{s,t}(x)$ in \eqref{TW} by any {\it regularized} flow $\theta^{(\eps)}_{s,t}(x)$
defined in \eqref{ODE0} below. When $\alpha\in(0,1)$, we choose the regularizing parameter  $\eps=(t-s)^{1/\alpha}$ since we need to use $\eps$ to compensate the
time singularity in the supercritical case.
For $\alpha\in(0,1]$, since $b$ is continuous in $x$, we can replace $\theta_{s,t}(x)$ in \eqref{TW} by any measurable {\it Peano} flow $\vartheta_{s,t}(x)$ of the ODE 
$\dot\vartheta_{s,t}(x)=b(t,\vartheta_{s,t}(x))$.
\er

\br\rm
When $b\equiv 0$ and $L^{(\alpha)}$ is a general $\alpha$-stable-like \textcolor{black}{generator of the form \eqref{GEN_STABLE_LIKE}}, it was proven in \cite{Li-So-Xi}
that the gradient estimate \eqref{GR0} holds for $\alpha\in(\frac12,2)$. See also \cite{Ch-Zh20} for more general 
\textcolor{black}{absolutely continuous L\'evy measures enjoying suitable scaling properties.
}
It seems that our gradient estimate \eqref{GR0} is the first result for SDEs of the form \eqref{SDE0} driven by a rotationally invariant $\alpha$-stable process
with $\alpha\in(0,\frac12]$.
\er

The paper is organized as follows. We give in Section \ref{SEC_PRELIM} some preliminary estimates needed for the main analysis. This concerns the mollified flow, some exit probabilities, convolution inequalities and the density of the \textit{proxy process} involved in the parametrix (which has a dynamic similar to \eqref{SDE0} with coefficient frozen along a suitable deterministic flow).
Section \ref{SEC_HK_SMOOTH_COEFF} is then {devoted} to the derivation of the two-sided bound and the fractional derivative estimate under {\bf (H$^\beta_b$)} and {\bf (H$^{\gamma}_a$)} when the coefficients are additionally supposed to be smooth and {\blue the drift bounded}. We specifically address the gradient estimate under those same assumptions in Section \ref{SEC_GRAD}. We eventually present in Section \ref{FINAL_SEC} some compactness arguments to derive the results of Theorem \ref{Main} under the sole conditions {\bf (H$^\beta_b$)}, {\bf (H$^{\gamma}_a$)}.
\section{Preliminaries}
\label{SEC_PRELIM}
\subsection{ODE flow}
We first present some basic properties about the solution $\theta_{s,t}(x)$ of the ODE \eqref{ODE1}. Since the drift coefficient therein 
depends on the initial time $s$, the following flow property does no longer hold:
$$
\theta_{r,t}\circ \theta_{s,r}(x)=\theta_{s,t}(x),\ \ s<r<t.
$$
However, the above flow property holds for the following regularized ODE:
\begin{align}\label{ODE0}
\dot\theta^{(\eps)}_{s,t}(x)=b_\eps(t,\theta^{(\eps)}_{s,t}(x)),\ \ \theta^{(\eps)}_{s,s}(x)=x,
\end{align}
for any fixed regularizing parameter $\eps>0 $.
Below we fix $\alpha\in(0,2)$ and always assume \eqref{BB0}.
The following lemma is easy.
\bl\label{Le41}
\begin{enumerate}[(i)]
\item For each $\eps>0$ and $s,t\geq 0$, $x\mapsto \theta^{(\eps)}_{s,t}(x)$ is a $C^1$-diffeomorphism and
\begin{align}\label{Inv}
(\theta^{(\eps)}_{s,t})^{-1}(y)=\theta^{(\eps)}_{t,s}(y). 
\end{align}
Moreover, for all $s,r,t\geq 0$, it holds that
\begin{align}\label{flow}
\theta^{(\eps)}_{s,t}(x)=\theta^{(\eps)}_{r,t}\circ\theta^{(\eps)}_{s,r}(x).
\end{align}

\item For all $\eps,\eps'>0$ and $s,t\geq 0$, $x\in\mR^d$, it holds that
\begin{align}\label{AQ4}
|\theta^{(\eps')}_{s,t}(x)-\theta^{(\eps)}_{s,t}(x)|\leq 2\kappa_0(\eps\vee\eps')^\beta
|t-s|\e^{\kappa_0\|\nabla\rho\|_{L^1}((\eps\vee\eps')^{\beta-1}+1)|t-s|},
\end{align}

\item For any $T>0$, there is a constant $C=C(T,d,\kappa_0)>0$ such that for all $s,t\in[0,T]$, $x,y\in\mR^d$ and $\eps=|t-s|^{1/\alpha}$,
\begin{align}\label{AQ3}
|\theta^{(\eps)}_{s,t}(x)-y|\asymp_C|x-\theta^{(\eps)}_{t,s}(y)|,\ \ |\theta^{(\eps)}_{s,t}(x)-\theta^{(\eps)}_{s,t}(y)|\asymp_C|x-y|.
\end{align}
\end{enumerate}
\el
\begin{proof}
(i) 
Note that by \eqref{ODE0}, for $0\le s<t $:
$$
\theta_{s,t}^{(\eps)}(x)=x+\int^t_s\textcolor{black}{b_{\eps}}(r,\theta_{s,r}^{(\eps)}(x))\dif r,\ \ 
\theta_{t,s}^{(\eps)}(y)=y-\int^t_s\textcolor{black}{b_{\eps}}(r,\theta_{t,r}^{(\eps)}(y))\dif r.
$$
Let $y=\theta_{s,t}^{(\eps)}(x)$. By the flow property, we have
\begin{align*}
y&=\big(\theta_{s,t}^{(\eps)}\big)^{-1}(y)+\int^t_s\textcolor{black}{b_{\eps}}(r,\theta_{s,r}^{(\eps)}\circ\big(\theta_{s,t}^{(\eps)}\big)^{-1}(y))\dif r\\
&=(\theta_{s,t}^{(\eps)})^{-1}(y)+\int^t_s\textcolor{black}{b_{\eps}}(r,\big(\theta_{r,t}^{(\eps)}\big)^{-1}(y))\dif r.
\end{align*}
Since the ODE has a unique solution, we immediately have $\big(\theta_{s,t}^{(\eps)}\big)^{-1}(y)=\theta_{t,s}^{(\eps)}(y)$. As for \eqref{flow}, it follows from \eqref{Inv}
and the flow property.
\\
\\
(ii) Without loss of generality, we assume $\eps'<\eps$. Since by \eqref{LL19} and \eqref{LL1},
\begin{align}\label{LL2}
|b_\eps(t,x)-b_{\eps'}(t,x)|\leq 2\kappa_0\eps^\beta,\ \ \|\nabla b_\eps\|_\infty\leq \kappa_0\|\nabla\rho\|_{L^1}(\eps^{\beta-1}+1),
\end{align}
by definition we have
\begin{align*}
|\theta^{(\eps')}_{s,t}(x)-\theta^{(\eps)}_{s,t}(x)|
&\leq\int^t_s |b_{\eps'}(r,\theta^{(\eps')}_{s,r}(x))-b_{\eps}(r,\theta^{(\eps')}_{s,r}(x))|\dif r
\\&\quad+\int^t_s |b_{\eps}(r,\theta^{(\eps')}_{s,r}(x)-b_{\eps}(r,\theta^{(\eps)}_{s,r}(x)|\dif r
\\&\leq 2\kappa_0\eps^\beta(t-s)+\kappa_0\|\nabla\rho\|_{L^1}(\eps^{\beta-1}+1)\int^t_s |\theta^{(\eps')}_{s,r}(x)-\theta^{(\eps)}_{s,r}(x)|\dif r.
\end{align*}
Using Gronwall's inequality, we obtain \eqref{AQ4}.
\\
\\
(iii) Without loss of generality, we assume $0\leq s<t\leq T$. 
Note that for $u\in[s,t]$,
$$
|\theta^{(\eps)}_{s,u}(x)-\theta^{(\eps)}_{s,u}(y)|\leq |x-y|+\int^u_s\|\nabla b_\eps(r,\cdot)\|_\infty|\theta^{(\eps)}_{s,r}(x)-\theta^{(\eps)}_{s,r}(y)|\dif r.
$$
For $\eps=|t-s|^{1/\alpha}$, it follows from the Gronwall inequality and \eqref{LL1} that
$$
|\theta^{(\eps)}_{s,t}(x)-\theta^{(\eps)}_{s,t}(y)|\leq \e^{\kappa_0\|\nabla\rho\|_{L^1}(|t-s|^{(\beta-1)/\alpha}+1)|t-s|}|x-y|.
$$
{\blue Similarly, from the backwards dynamics, one can derive that for \textcolor{black}{$x',y'\in \R^d $},
$$
|\theta^{(\eps)}_{t,s}(\textcolor{black}{x'})-\theta^{(\eps)}_{t,s}(\textcolor{black}{y'})|\leq \e^{\kappa_0\|\nabla\rho\|_{L^1}(|t-s|^{(\beta-1)/\alpha}+1)|t-s|}|\textcolor{black}{x'-y'}|.
$$
Thus, by \eqref{Inv} \textcolor{black}{(inverse flow property)} and putting $x'=\theta^{(\eps)}_{\textcolor{black}{s,t}}(x)$ and $y'=\theta^{(\eps)}_{\textcolor{black}{s,t}}(y)$ we obtain the second equivalence in \eqref{AQ3}.
\textcolor{black}{The first equivalence follows from the second one  replacing therein $y $ by $\theta_{t,s}^{\eps}(y) $ and using \eqref{Inv}}. 
}
\end{proof}
The following result is a consequence of the above lemma, which plays a crucial role below.
\bl\label{Le22}
(i) For each $s,t\geq 0$, the map $x\mapsto \theta_{s,t}(x)$, where $(\theta_{s,u}(x))_{u\in [s,t]}$ solves \eqref{ODE1}, 
is a $C^1$-diffeomorphism and there is a constant $C_0=C_0(T,\Theta)>0$
such that
$$
|\det(\nabla\theta^{-1}_{s,t}(x))-1|\le C_0|t-s|^{(\alpha+\beta-1)/\alpha}.
$$
(ii) For any $T>0$, there is a constant $C_1=C_1(T,\Theta)\geq 1$
such that for all $s,t\in[0,T]$ and $x,y\in\mR^d$,
\begin{align}\label{KH2}
|t-s|^{1/\alpha}+|\theta_{s,t}(x)-y|\asymp_{C_1}|t-s|^{1/\alpha}+|x-\theta_{t,s}(y)|.
\end{align}
(iii) For any $T>0$, there is a constant $C_2=C_2(T,\Theta)>0$ such that for all $s,r,t\in[0,T]$ and $x\in\mR^d$,
\begin{align}\label{KH1}
|\theta_{s,t}(x)-\theta_{r,t}\circ\theta_{s,r}(x)|\le C_2 |(r\vee s\vee t)-(r\wedge s\wedge t)|^{1/\alpha}.
\end{align}
\el
\begin{proof}
(i) It is well known that {\green (cf. \cite[Proposition 1.2]{MB02})}
$$
\det(\nabla\theta_{s,t}(x))=1+\int^t_s\div b_{|r-s|^{1/\alpha}}(r,\theta_{s,r}(x))\det(\nabla\theta_{s,r}(x))\dif r.
$$
Thus
$$
\det(\nabla\theta_{s,t}(x))=\exp\left\{\int^t_s\div b_{|r-s|^{1/\alpha}}(r,\theta_{s,r}(x))\dif r\right\}.
$$
The desired estimate follows by \eqref{LL01} and $\nabla\theta^{-1}_{s,t}(x)=(\nabla\theta_{s,t})^{-1}(\theta^{-1}_{s,t}(x))$.
\\
\\
(ii) Fix $s<t$. For $u\in[s,t]$, by definition we have
$$
\theta_{s,u}(x)=x+\int^u_sb_{|r-s|^{1/\alpha}}(r,\theta_{s,r}(x))\dif r,
$$
and for $\eps=|t-s|^{1/\alpha}$,
$$
\theta^{(\eps)}_{s,u}(x)=x+\int^u_sb_{\eps}(r,\theta^{(\eps)}_{s,r}(x))\dif r.
$$
By \eqref{LL2} with $\eps'=|r-s|^{1/\alpha}$ and $\eps=|t-s|^{1/\alpha}$, we have for all $u\ge s $,
\begin{align*}
|\theta_{s,u}(x)-\theta^{(\eps)}_{s,u}(x)|
&\le \int^u_s\Big|b_{|r-s|^{1/\alpha}}(r,\theta_{s,r}(x))-b_{\eps}(r,\theta^{(\eps)}_{s,r}(x))\Big|\dif r
\\&\lesssim |t-s|^{\beta/\alpha+1}+(t-s)^{(\beta-1)/\alpha}\int^u_s|\theta_{s,r}(x)-\theta^{(\eps)}_{s,r}(x)|\dif r,
\end{align*}
which yields by Gronwall's inequality that
\begin{align}\label{AQ5}
|\theta_{s,t}(x)-\theta^{(\eps)}_{s,t}(x)|\lesssim\e^{(t-s)^{(\beta-1)/\alpha+1}} |t-s|^{\beta/\alpha+1}\lesssim |t-s|^{1/\alpha},
\end{align}
where the second inequality is due to $\alpha+\beta>1$.
Thus, by \eqref{AQ3}  and \eqref{AQ5}, we have
\begin{align*}
|\theta_{s,t}(x)-y|&\leq|\theta^{(\eps)}_{s,t}(x)-y|+|\theta_{s,t}(x)-\theta^{(\eps)}_{s,t}(x)|
\\&\lesssim |x-\theta^{(\eps)}_{s,t}(y)|+|t-s|^{1/\alpha}\lesssim |x-\theta_{t,s}(y)|+|t-s|^{1/\alpha}.
\end{align*}
The right hand side inequality of \eqref{KH2} follows. By symmetry, we also have the left hand side inequality.
\\
\\
(iii) Let $s,r,t\geq 0$ and $\eps:=|r\vee s\vee t-r\wedge s\wedge t|^{1/\alpha}$. By \eqref{flow} we have
\begin{align*}
|\theta_{s,t}(x)-\theta_{r,t}\circ\theta_{s,r}(x)|
&\leq|\theta_{s,t}(x)-\theta^{(\eps)}_{s,t}(x)|
+|\theta^{(\eps)}_{r,t}\circ\theta^{(\eps)}_{s,r}(x)-\theta^{(\eps)}_{r,t}\circ\theta_{s,r}(x)|
\\&+|\theta^{(\eps)}_{r,t}\circ\theta_{s,r}(x)-\theta_{r,t}\circ\theta_{s,r}(x)|.
\end{align*}
The desired estimate \eqref{KH1} again follows by \eqref{AQ3} and \eqref{AQ5}.
\end{proof}
\br\rm
By \eqref{KH1}, we have
\begin{align}\label{KH11}
|x-\theta_{t,s}\circ\theta_{s,t}(x)|\lesssim_{C_2} |t-s|^{1/\alpha},
\end{align}
and by \eqref{KH2},
\begin{align}\label{KH22}
|t-s|^{1/\alpha}+|x-y|\asymp_{C_1}|t-s|^{1/\alpha}+|\theta_{s,t}(x)-\theta_{t,s}(y)|.
\end{align}
Put it differently, we have an \textit{approximate} flow property for the ODE \eqref{ODE1}. Namely, the flow property holds up to an additive time factor which has the same magnitude as the current typical time (self-similarity index of the driving process).
\er
\subsection{Probability estimates}

We need the following master formula.
\bl\label{Le24}
(L\'evy system) Let $X_t:=X_{0,t}$ be any solution of SDE \eqref{SDE0}.
For any nonnegative measurable function $f: \mR_+\times\mR^d\times\mR^d\to\mR_+$ and finite stopping time $\tau$, 
$$
\mE\sum_{r\in(0,\tau]}f(r, X_{r-}, \Delta X_r)=\mE\int^\tau_0\!\!\!\int_{\mR^d}f(r, X_{r-}, z)\frac{\kappa(r,X_{r-},z)}{|z|^{d+\alpha}}\dif z\dif r,
$$
where $\Delta X_r:=X_r-X_{r-}$ and $\kappa(r,x,z)$ is defined in \eqref{KA1}.
\el
\textcolor{black}{Let us emphasize that the L\'evy system is well known for a time-homogeneous process, we can e.g. refer to Chapter 2.4 in \cite{bott:schi:wang:13} and references therein. The extension to the current framework under the uniform ellipticity condition  is rather natural. We provide a proof for the sake of completeness.}
\begin{proof}
Let $N(\dif t,\dif z)$ be the counting measure associated with $L^{(\alpha)}_t$, i.e.,
$$
N((0,t]\times E):=\sum_{s\in[0,t]}\1_E(\Delta L^{(\alpha)}_s),\ { E\in\sB(\mR^d)}.
$$
Noting that 
$$
\Delta X_t=a(t,X_{t-})\Delta L^{(\alpha)}_t,
$$
we have for any $\eps>0$,
$$
\sum_{r\in(0,\tau]}f(r, X_{r-}, \Delta X_r)
\1_{ |\Delta X_r|\geq\eps}
=\int^\tau_0\!\!\!
\int_{|{ a(r,X_{r-})z}|\geq\eps}
f(r, X_{r-}, a(r,X_{r-})z)
N(\dif r,\dif z).
$$
Since the compensated measure of $N(\dif t,\dif z)$ is $\frac{\dif z\dif t}{|z|^{d+\alpha}}$, 
by the change of variable, we have
\begin{align*}
\mE\sum_{r\in(0,\tau]}f(r, X_{r-}, \Delta X_r)
\1_{|\Delta X_r|\geq\eps}
&=\mE\int^\tau_0\!\!\!
\int_{|{ a(r,X_{r-})z}|\geq\eps}
f(r, X_{r-}, a(r,X_{r-})z)\frac{\dif z\dif r}{|z|^{d+\alpha}}
\\&=\mE\int^\tau_0\!\!\!\int_{|z|\geq\eps}f(r, X_{r-}, z)\frac{\kappa(r,X_{r-},z)}{|z|^{d+\alpha}}\dif z\dif r,
\end{align*}
where $\kappa(r,x,z)$ is given in \eqref{KA1},
which in turn gives the desired formula by the monotone convergence theorem.
\end{proof}

Fix $(s,x)\in\mR_+\times\mR^d$. For $\eta>0$, define the stopping time
\begin{align}\label{DH1}
\tau^\eta_{s,x}:=\inf\Big\{t>s: |X_{s,t}(x)-\theta_{s,t}(x)|>\eta\Big\},
\end{align}
which corresponds to the exit time of the diffusion from a tube around the deterministic ODE introduced in \eqref{ODE1}.
We now give a \textit{tube estimate} which roughly says that, for a given spatial threshold $\eta$, the probability that the difference between the process $X_{s,\cdot}(x)$ and the deterministic regularized flow $\theta_{s,\cdot}(x) $ leaves the tube of radius $\eta $ before a certain fraction $\varepsilon \eta^\alpha $ of the corresponding typical time scale $\eta^\alpha $ is somehow \textit{small}.
\bl\label{Le23}
Under {\bf (H$^\beta_b$)} and {\bf (H$^0_a$)},  for any $T>0$,
there is an $\eps\in(0,1)$ only depending on $T,\Theta$ such that
for all $(s,x)\in\mR_+\times\mR^d$ and $\eta\in(0,T^{1/\alpha}]$, 
$$
\mP(\tau^\eta_{s,x}<s+\eps\eta^\alpha)\leq 1/2.
$$
\el
\begin{proof}
Without loss of generality, we assume $(s,x)=(0,0)$, and for simplicity write
$$
X_t:=X_{0,t}(0),\ \ \theta_t:=\theta_{0,t}(0),\ \ \tau:=\tau^\eta_{0,0}.
$$
Let $f\in C^2_b(\mR^d)$ with $f(0)=0$ and $f(x)=1$ for $|x|>1$. For $\eta>0$, set
$$
f_\eta(x):=f(x/\eta),\ \ u:=\eps\eta^\alpha.
$$
Note that
$$
Y_t:=X_t-\theta_t=\int^t_0(b(r,X_r)-b_{r^{1/\alpha}}(r,\theta_r))\dif r+\int^t_0a(r,X_{r-})\dif L^{(\alpha)}_r.
$$
By It\^o's formula, we have
\begin{align*}
\mE f_\eta(Y_{u\wedge\tau})&=\mE\int^{u\wedge\tau}_0\!\!\Big[(b(r,X_r)-b_{r^{1/\alpha}}(r,\theta_r))\cdot\nabla f_\eta(Y_r)
+\tfrac12\cL_r f_\eta(Y_r)\Big]\dif r.
\end{align*}
Note that
\begin{align*}
|b(r,X_r)-b_{r^{1/\alpha}}(r,\theta_r)|
&\leq|b(r,X_r)-b(r,\theta_r)|+|b(r,\theta_r)-b_{r^{1/\alpha}}(r,\theta_r)|
\\&\underset{\eqref{BB0},\eqref{LL2}}{\leq}\kappa_0(|X_r-\theta_r|^\beta+|X_r-\theta_r|)+\kappa_0r^{\beta/\alpha}.
\end{align*}
Also, with the notation of \eqref{DEF_DELTA_2},
$$
|\delta^{(2)}_{f_\eta}(x;z)|\leq (|z|^2\|\nabla^2 f_\eta\|_\infty)\wedge(4\|f_\eta\|_\infty)\lesssim (|z|^2/\eta^2)\wedge 1.
$$
Hence,
\begin{align*}
\mE f_\eta(Y_{u\wedge\tau})&\lesssim\mE\int^{u}_0\left[(|Y_r|^\beta+|Y_r|+r^{\beta/\alpha})\cdot\frac{\1_{|Y_r|\leq\eta}}{\eta}
+\int_{\mR^d}\frac{(|z|^2/\eta^2)\wedge 1}{|z|^{d+\alpha}}\dif z\right]\dif r
\\&\lesssim\int^{u}_0(\eta^{\beta-1}+1+r^{\beta/\alpha}\eta^{-1}+\eta^{-\alpha})\dif r
\\&\underset{u=\eps\eta^\alpha}{\lesssim} \eps(\eta^{\alpha+\beta-1}+\eta^\alpha+ \eps^{\beta/\alpha}\eta^{\alpha+\beta-1}+1)
\lesssim \eps( T^{1+(\beta-1)/\alpha}+T+1),
\end{align*}
where the last inequality is due to $\alpha+\beta>1$ and $\eta\leq T^{1/\alpha}$.  Importantly, the implicit constant is
independent of $\eps$. Note that
\begin{align*}
\mP(\tau<u)=\mE\1_{\tau<u}\leq\mE f_{\eta}\big(Y_{u\wedge\tau}\big)\lesssim \eps(T^{1+(\beta-1)/\alpha}+T+1).
\end{align*}
The desired estimate follows by choosing $\eps$ small enough.
\end{proof}
The following lemma will be used to show the lower bound estimate of the heat kernel. 
It gives a lower bound estimate for the probability that, considered two points $x,y$ in the \textit{off-diagonal regime} between times $s$ and $t$, namely such that $ |x-\theta_{t,s}(y)|\ge K(t-s)^{1/\alpha}$, after a time $\varepsilon (t-s)$ with $\eps $ as in Lemma  \ref{Le23}, the stochastic forward transport of $x$ by the SDE, i.e. $X_{s,s+\eps(t-s)}(x) $ and the backward deterministic transport of $y$ by the regularized flow $\theta_{t,s+\eps(t-s)}(y) $ belong to a \textit{diagonal} tube, with radius $K(t-s)^{1/\alpha} $,
where again  $(t-s)^{1/\alpha} $ corresponds to the current typical scale between times $s$ and $t$.

\bl\label{Le26}
Suppose that {\bf (H$^\beta_b$)} and {\bf (H$^0_a$)} hold. Let $\eps\in(0,1)$ be as in Lemma \ref{Le23}. For any $T>0$,
there are constants $c_0\in(0,1),K\geq 1$ depending only on $T,\Theta$ 
such that for all $0\leq s<t\leq T$ and $|x-\theta_{t,s}(y)|\geq K(t-s)^{1/\alpha}$,
$$
\mP\Big(|X_{s,s+\eps(t-s)}(x)-\theta_{t,s+\eps(t-s)}(y)|\leq K(t-s)^{1/\alpha}\Big)\geq \frac{c_0(t-s)^{1+d/\alpha}}{|x-\theta_{t,s}(y)|^{d+\alpha}}.
$$
\el
This Lemma will be crucial for the lower bound estimate of the heat kernel, since it precisely gives the control needed for a chaining argument, see Theorem \ref{T8} below. As opposed to the continuous case, for SDEs driven by pure jump processes, a single intermediate time, associated with a \textit{large jump}, is needed for the chaining. Roughly speaking between times $s+\eps(t-s)$ and $t$ we will use the global \textit{diagonal} bound of order $(t-s)^{-d/\alpha} $, since $ \eps$ is meant to be \textit{small enough}, and the above Lemma controls the probability that the process enters a good neigborhood of the backward flow to do so. The lower bound is the sought one in the sense that when 
{ multiplying it} by $(t-s)^{-d/\alpha}$ exactly makes the expression in \eqref{TW}, 
$ (t-s)((t-s)^{1/\alpha}+|\theta_{s,t}(x)-y|)^{-(d+\alpha)}\asymp_C (t-s)|\theta_{s,t}(x)-y|^{-(d+\alpha)}$ appear since $|x-\theta_{s,t}(y)|\ge K(t-s)^{1/\alpha} $ (\textit{off-diagonal regime}). 

\begin{proof}[Proof of Lemma \ref{Le26}]
Without loss of generality, we assume $s=0$ and for simplicity, we write
$$
\eta:=t^{1/\alpha},\ u:=\eps\eta^\alpha=\eps t,\ \ X_r(x):=X_{0,r}(x).
$$
Define a stopping time
$$
\sigma:=\inf\Big\{r>0: |X_{r}(x)-\theta_{t,r}(y)|\leq \eta\Big\}.
$$
By the right continuity of $r\mapsto X_{r}(x)-\theta_{t,r}(y)$, one sees that
$$
|X_{\sigma}(x)-\theta_{t,\sigma}(y)|\leq \eta, \  a.s.
$$
In particular, for $\sigma\leq u$, by \eqref{KH22} and \eqref{KH1}, there is a constant $C_0=C_0(\Theta)\geq 1$ such that
\begin{align*}
|X_{u}(x)-\theta_{t,u}(y)|
&\leq |X_{\sigma,u}(X_{\sigma}(x))-\theta_{\sigma,u}(X_{\sigma}(x))|
\\&\quad+|\theta_{\sigma,u}(X_{\sigma}(x))-\theta_{\sigma,u}(\theta_{t,\sigma}(y))|
\\&\quad+|\theta_{\sigma,u}(\theta_{t,\sigma}(y))-\theta_{t,u}(y)|
\\&\leq |X_{\sigma,u}(X_{\sigma}(x))-\theta_{\sigma,u}(X_{\sigma}(x))|+C_0\eta.
\end{align*}
Let $K\geq C_0+1$. Then
$$
\big\{|X_{\sigma,u}(X_{\sigma}(x))-\theta_{\sigma,u}(X_{\sigma}(x))|<\eta\big\}
\subset\big\{|X_{u}(x)-\theta_{t,u}(y)|\leq K\eta\big\}.
$$
Thus, by the strong Markov property, we have
\begin{align*}
\mP(|X_{u}(x)-\theta_{t,u}(y)|\leq K\eta)
&\geq \mP(\sigma\leq u; |X_{u}(x)-\theta_{t,u}(y)|\leq K\eta)
\\&\geq \mP\left(\sigma\leq u; |X_{\sigma,u}(X_{\sigma}(x))-\theta_{\sigma,u}(X_{\sigma}(x))|<\eta\right)
\\&\geq \mP(\sigma\leq u)\inf_{(s,z)\in[0,u]\times\mR^d}\mP\left(|X_{s,u}(z)-\theta_{s,u}(z)|<\eta\right).
\end{align*}
Let $\tau^\eta_{s,z}$ be defined by \eqref{DH1}. 
By Lemma \ref{Le23} we have
\begin{align}\label{Tau}
\mP\left(|X_{s,u}(z)-\theta_{s,u}(z)|\geq \eta\right)\leq \mP(\tau^\eta_{s,z}\leq u)
\leq \mP(\tau^\eta_{s,z}\leq s+\eps\eta^\alpha)\leq 1/2,
\end{align}
which implies that
\begin{align*}
\inf_{(s,z)\in[0,u]\times\mR^d}\mP\left(|X_{s,u}(z)-\theta_{s,u}(z)|<\eta\right)\geq 1/2
\end{align*}
and
\begin{align}\label{AA9}
\mP(|X_{u}(x)-\theta_{t,u}(y)|\leq K\eta)\geq \mP(\sigma\leq u)/2.
\end{align}
Next we need to obtain a lower bound estimate for $\mP(\sigma\leq u)$. Let $\tau:=\tau^\eta_{0,x}$.
For $r<u\wedge\tau$, by \eqref{KH2}, 
 there are constants $c_0, C_1>0$ such that
\begin{align*}
|X_{r}(x)-\theta_{t,r}(y)|&\geq|\theta_{0,r}(x)-\theta_{t,r}(y)|-|X_{r}(x)-\theta_{0,r}(x)|\no
\\&\geq c_0|x-\theta_{t,0}(y)|-C_1 t^{1/\alpha}-\eta.
\end{align*}
In particular, if we choose $K\geq (C_0+1)\vee ((C_1+2)/c_0)$, then since  by assumption
$$
|x-\theta_{t,0}(y)|\geq K\eta,
$$
it holds that for $r<u\wedge\tau$,
\begin{align}
|X_{r}(x)-\theta_{t,r}(y)|\geq c_0K\eta-(C_1+1)\eta>\eta.\label{DG2}
\end{align}
Thus we have
$$
\1_{\{|X_{u\wedge\tau}(x)-\theta_{t,u\wedge\tau}(y)|\leq\eta\}}
=\sum_{r\in(0,u\wedge\tau]}\1_{\{|X_r(x)-\theta_{t,r}(y)|\leq\eta\}},
$$
i.e. we have at most one term in the above summand. We then derive from Lemma \ref{Le24},
\begin{align*}
\mP{\{|X_{u\wedge\tau}(x)-\theta_{t,u\wedge\tau}(y)|\leq\eta\}}
&=\mE\sum_{r\in(0,u\wedge\tau]}\1_{\{|X_r(x)-\theta_{t,r}(y)|\leq\eta\}}\\
&=\mE\int^{u\wedge\tau}_0\!\!\!\int_{|z-\theta_{t,r}(y)|\leq\eta}\frac{\kappa(r, X_{r-}(x), z-X_{r-}(x))}{|z-X_{r-}(x)|^{d+\alpha}}\dif z\dif r.
\end{align*}
On the other hand, noting that for $r<u\wedge\tau$ and $|z-\theta_{t,r}(y)|\leq\eta$, we get
\begin{align*}
|z-X_{r-}(x)|&\leq|z-\theta_{t,r}(y)|+|\theta_{t,r}(y)-\theta_{0,r}(x)|+|\theta_{0,r}(x)-X_{r-}(x)|
\\&\leq \eta+C_2|x-\theta_{t,0}(y)|+C_2 t^{1/\alpha}+\eta\leq C_2|x-\theta_{t,0}(y)|+C_3\eta,
\end{align*}
using as well \eqref{KH2} for the last but one inequality. Since
$$
\{|X_{u\wedge\tau}(x)-\theta_{t,u\wedge\tau}(y)|\leq\eta\}\subset \{\sigma\leq u\},
$$
we further have
\begin{align}
\mP\{\sigma\leq u\}&\geq \mP{\{|X_{u\wedge\tau}(x)-\theta_{t,u\wedge\tau}(y)|\leq\eta\}}\no\\
&\geq\mE\int^{u\wedge\tau}_0\int_{|z-\theta_{t,r}(y)|\leq\eta}\frac{\kappa_1^{-1}}
{(C_2|x-\theta_{t,0}(y)|+C_3\eta)^{d+\alpha}}\dif z\dif r\no
\\&=\frac{\mE(u\wedge\tau)\kappa_1^{-1}\eta^d\cdot {\rm Vol}(B_1)}{(C_2|x-\theta_{t,0}(y)|+C_3\eta)^{d+\alpha}}
\geq\frac{c_0 t^{1+d/\alpha}}{|x-\theta_{t,0}(y)|^{d+\alpha}},\label{AA10}
\end{align}
where the last step is due to $|x-\theta_{t,0}(y)|\geq K\eta$ and 
$$
\mE (u\wedge\tau)\geq u\mP(\tau>u)\stackrel{\eqref{Tau}}{\geq} u/2=\eps t/2.
$$
Combining \eqref{AA9} and \eqref{AA10}, we obtain the desired estimate.
\end{proof}
\subsection{Convolution inequalities}
This Section is dedicated to some useful convolution controls associated with functions that are known to be upper-bounds of the isotropic stable density and its gradient, see e.g. \cite{kolo:00}, \cite{Bo-Ja}. Though a bit technical, these results will turn out to be crucial in order to control the parametrix series representation of the density and its gradient (see e.g.  equation \eqref{Eq1}, Lemma \ref{Le31} and Theorem \ref{T8} below).

For $\eta\in(0,2)$ and $(t,x)\in\mR_+\times\mR^d$, let
$$
\varrho^{(\eta)}(x):=(1+|x|)^{-d-\eta},\ \ 
\varrho^{(\eta)}(t,x):=t^{-d/\alpha}\varrho^{(\eta)}(t^{-1/\alpha}x).
$$
For $\beta\geq 0$ and $\gamma\in\mR$, we introduce the following functions for later use
\begin{align}\label{DEF_VAR_RHO}
\varrho^{(\eta)}_{\beta,\gamma}(t,x)&:=(1\wedge (t^{1/\alpha}+|x|))^\beta t^{(\gamma-\eta)/\alpha}\varrho^{(\eta)}(t,x)
\end{align}
and  for $0\leq s<t$,
\begin{align}
\label{DEF_PHI}
\phi^{(\eta)}_{\beta,\gamma}(s,x,t,y):=\varrho^{(\eta)}_{\beta,\gamma}(t-s,x-\theta_{t,s}(y)).
\end{align}
Note that
$$
\varrho^{(\eta)}_{\beta,\gamma}(t,x)=\frac{(1\wedge (t^{1/\alpha}+|x|))^\beta t^{\gamma/\alpha}}{(t^{1/\alpha}+|x|)^{d+\eta}}.
$$
For $T>0$, by \eqref{KH2} we have for $(s,x,t,y)\in\mD_T$,
\begin{align}\label{ES44}
\phi^{(\eta)}_{\beta,\gamma}(s,x,t,y)\asymp \varrho^{(\eta)}_{\beta,\gamma}(t-s,\theta_{s,t}(x)-y),
\end{align}
and for $\beta\in[0,\eta]$,
\begin{align}
\int_{\mR^d}\phi^{(\eta)}_{\beta,\gamma}(s,x,t,y)\dif y
&\lesssim\int_{\mR^d}\varrho^{(\eta)}_{\beta,\gamma}(t-s,y)\dif y
\lesssim  (t-s)^{\frac{\beta+\gamma-\eta}{\alpha}}.\label{ES4}
\end{align}
\textcolor{black}{Observe as well importantly that, from \eqref{DEF_VAR_RHO},  for $0\le \beta'\le \beta $, $0\le \gamma'\le \gamma $
and  $(s,x,t,y)\in\mD_T$,
\begin{align}
\label{DOM_STAB}
\phi^{(\eta)}_{\beta,\gamma}(s,x,t,y)\lesssim \phi^{(\eta)}_{\beta',\gamma'}(s,x,t,y).
\end{align}
}
For two functions $f,g$ on $\mD_\infty$, we write
$$
(f\odot g)_r(s,x,t,y):=\int_{\mR^d} f(s,x,r,z)g(r,z,t,y)\dif z
$$
and
$$
(f\otimes g)(s,x,t,y):=\int^t_s(f\odot g)_r(s,x,t,y)\dif r.
$$
The following lemma is the same as in \cite[Lemma 2.1]{Ch-Zh}.
\bl\label{Le09}  Fix $\alpha\in (0,2) $.
For any $\beta_1,\beta_2\in[0,\frac{\alpha}{4}]$ and $T>0$,
there is a  $C=C(T,\Theta,\beta_1,\beta_2)>0$ such that
for all $\gamma_1>-\beta_1$ and $\gamma_2>-\beta_2$, $r\in[s,t]$ and $x,y\in \R^d $,
\begin{align}\label{eq30}
\big(\phi^{(\alpha)}_{\beta_1,0}\odot\phi^{(\alpha)}_{\beta_2,0}\big)_r(s,x,t,y)
\lesssim_{C} \big((r-s)^{\frac{\beta_1-\alpha}{\alpha}}+(t-r)^{\frac{\beta_2-\alpha}{\alpha}}\big)
\phi^{(\alpha)}_{\beta_1\wedge\beta_2,0}(s,x,t,y)
\end{align}
and
\begin{align}\label{eq31}
\phi^{(\alpha)}_{\beta_1,\gamma_1}\otimes\phi^{(\alpha)}_{\beta_2,\gamma_2} (s,x,t,y)
\lesssim_{C} \cB(\tfrac{\beta_1+\gamma_1}{\alpha},\tfrac{\beta_2+\gamma_2}{\alpha})
\phi^{(\alpha)}_{\beta_1\wedge\beta_2,\beta_1+\beta_2+\gamma_1+\gamma_2}(s,x,t,y), 
\end{align}
where $\cB(\gamma,\beta)$ is the usual Beta function defined by
$$
\cB(\gamma,\beta):=\int^1_0(1-s)^{\gamma-1}s^{\beta-1}\dif s,\ \ \gamma,\beta>0.
$$
\el
\begin{proof}
We follow the proof in \cite{Ch-Zh20}.
Let $\ell(u):=\frac{u^{d+\alpha}}{1\wedge u^{\beta}}$. It is easy to see that,  as soon as $d+\alpha>\beta $, $\ell$ is increasing on $\mR_+$ and
for any $\lambda\geq 1$,
\begin{align}\label{KP10}
\ell(\lambda u)\leq \lambda^{d+\alpha}\ell(u).
\end{align}
Hence,
\begin{align}\label{KP9}
\ell(u+w)\leq \ell(2(u\vee w))\leq 2^{d+\alpha}\ell(u\vee w)\leq 2^{d+\alpha}(\ell(u)+\ell(w)).
\end{align}
Now for $r\in[s,t]$ and $x,y\in\mR^d$, since 
$$
|t+s|^{1/\alpha}+|x+y|\leq 2^{1/\alpha}\big(|s|^{1/\alpha}+|x|+|t|^{1/\alpha}+|y|\big),
$$
by \eqref{KP10} and \eqref{KP9}, we have
\begin{align*}
\ell(|t+s|^{1/\alpha}+|x+y|)\lesssim_C\ell(|s|^{1/\alpha}+|x|)+\ell(|t|^{1/\alpha}+|y|).
\end{align*}
In particular,
\begin{align*}
\frac{((t+s)^{1/\alpha}+|x+y|)^{d+\alpha}}{1\wedge ((t+s)^{1/\alpha}+|x+y|)^{\beta_1\wedge\beta_2}}
\lesssim\frac{(s^{1/\alpha}+|x|)^{d+\alpha}}{1\wedge (s^{1/\alpha}+|x|)^{\beta_1}}
+\frac{(t^{1/\alpha}+|y|)^{d+\alpha}}{1\wedge (t^{1/\alpha}+|y|)^{\beta_2}}.
\end{align*}
Hence,
\begin{align*}
&\frac{1\wedge (s^{1/\alpha}+|x|)^{\beta_1}}{(s^{1/\alpha}+|x|)^{d+\alpha}}\times\frac{1\wedge (t^{1/\alpha}+|y|)^{\beta_2}}{(t^{1/\alpha}+|y|)^{d+\alpha}}\\
&\lesssim_C\left[\frac{1\wedge (s^{1/\alpha}+|x|)^{\beta_1}}{(s^{1/\alpha}+|x|)^{d+\alpha}}
+\frac{1\wedge (t^{1/\alpha}+|y|)^{\beta_2}}{(t^{1/\alpha}+|y|)^{d+\alpha}}\right]\\
&\quad\times\frac{1\wedge ((t+s)^{1/\alpha}+|x+y|)^{\beta_1\wedge\beta_2}}{((t+s)^{1/\alpha}+|x+y|)^{d+\alpha}}.
\end{align*}
By this, the desired estimates follow by \eqref{ES44}, \eqref{ES4} and \eqref{KH1}.
\end{proof}

\subsection{Density estimate } 
\label{SEC_DENS_EST}
Let $a:\mR_+\to\mR^d\otimes\mR^d$ be a measurable $d\times d$-matrix-valued function satisfying 
the non-degeneracy condition
\begin{align}\label{AA4}
\kappa_1^{-1}|\xi|^2\leq |a(s)\xi|^2\leq \kappa_1|\xi|^2.
\end{align}
Fix $\alpha\in (0,2)$ and consider the following jump process
\begin{equation}\label{SUB_INTEGRAL}
X^a_{s,t}:=\int^t_sa(r)\dif W_{S_r},
\end{equation}
where $ W$ is a $d$-dimensional Brownian motion and $S$ is an $\alpha/2$-stable subordinator independent from $W$, both defined on some probability space $(\Omega,\mathcal F,\mathbb P) $. Note that
$$
X^a_{s,t}\stackrel{(d)}{=}(t-s)^{1/\alpha}X^{\tilde a}_{0,1},
$$
where
$$
\tilde a(r):=a(s+r(t-s)).
$$

We have the following lemma that can be derived from the approach initially used in \cite{Be} (see also 
 \cite{Bo-Ja}).
 We provide below a proof for completeness.
\bl\label{Le27}
For any $0\leq s<t<\infty$, $X^a_{s,t}$ has a smooth density $p_{s,t}^{a}(x)$ with the scaling property 
\begin{align}\label{SC0}
p_{s,t}^{a}(x)=(t-s)^{-d/\alpha} p_{0,1}^{\tilde a}((t-s)^{-1/\alpha}x), 
\end{align}
which enjoys the following estimates:
\begin{align}
p_{s,t}^{a}(x)\asymp_{C_0}\varrho^{(\alpha)}_{0,\alpha}(t-s,x),\label{Es0}
\end{align}
and for any $j\in\mN$,
\begin{align}
|\nabla^j p_{s,t}^{a}(x)|&\lesssim_{C_j}\varrho^{(\alpha+j)}_{0,\alpha}(t-s,x),\label{Es2}
\end{align}
where the constants  only depend on $\kappa_1,d,\alpha$. Moreover, suppose that the integrand in \eqref{SUB_INTEGRAL} writes as $a_\xi(r)$ and smoothly depends on the parameter $\xi\in\mR^d$ so that \eqref{AA4} holds uniformly and
$\sup_{r,\xi}|\nabla^k_\xi a_\xi(r)|<\infty$ for any $k\in\mN$. 
Let $p_{s,t}^{a_\xi}$ be the density of the integral in \eqref{SUB_INTEGRAL} associated with $a_\xi$. Then we have
for $k\in\mN$ and $j\in\mN_0$,
\begin{align}
|\nabla^k_\xi\nabla^j_x p_{s,t}^{a_\xi}(x)|\lesssim_{C_{j,k}} \varrho^{(\alpha+j)}_{0,\alpha}(t-s,x).\label{Es22}
\end{align}
Importantly, this last bound means that, { the} differentiation w.r.t. the parameter $\xi $ appearing in the  \textit{diffusion} coefficient $a_\xi$ does not yield an additional time singularity.
 \el
\begin{proof}
The two sided estimate \eqref{Es0} is well known (see e.g. \cite{Ch-Zh1}). 
We show \eqref{Es2}.
Without loss of generality, we assume $s=0$ and write
$$
X_{t}:=\int^t_0a(r)\dif W_{S_r}.
$$
Fix a c\`adl\`ag path $\ell_s$. Consider the following Gaussian random variable:
$$
X^\ell_{t}:=\int^t_0a(r)\dif W_{\ell_r}.
$$
It has a density 
\begin{equation}\label{UNDERLYING_SUB_GAUSS}
 g^{a,\ell}_{t}(x)=(2\pi)^{-d/2}\sqrt{\det\big( (\cC_{t}^{a,\ell})^{-1}\big)}\exp\{-\<\big(\cC_t^{a,\ell}\big)^{-1}x,x\>/2\},
\end{equation}
where
$$
\cC_t^{a,\ell}:=\int^t_0(aa^*)(r)\dif \ell_r.
$$
From the non-degeneracy assumption \eqref{AA4}, we have
$$
\<\big(\cC_t^{a,\ell}\big)^{-1}x,x\>\asymp|x|^2/\ell_t,\ \ \det\big( (\cC_t^{a,\ell})^{-1}\big)\asymp \ell_t^{-d},
$$
and
$$
|\nabla  g^{a,\ell}_{t}(x)|\lesssim |x|/\ell_t\exp\{-\lambda|x|^2/\ell_t\}.
$$
The density $p_{0,t}^{a}(x)=:p_{t}^{a}(x)$ of $X_{t}$ is given by
\begin{equation}\label{SUB_STABLE_DENS}
p_{t}^{a}(x)=\mE  g^{a,S}_{t}(x).
\end{equation}
The bound of \eqref{Es2} is direct from the Fourier representation of the density when $|x|\le t^{1/\alpha}$. 
On the other hand, for $|x|>t^{1/\alpha} $,  from the global bound on the law  of the subordinator
$$
\mu_{S_t}(\dif r):=\mP\circ S_t^{-1}(\dif r)\lesssim t\, r^{-\alpha/2-1}\dif r, 
$$
it readily follows that
\begin{align*}
|\nabla  p_{t}^a(x)|\leq\mE |\nabla  g^{a,S}_{t}(x)|
&\lesssim|x|\mE(S_t^{-d/2-1}\exp\{-\lambda|x|^2/S_t\})<+\infty.
\end{align*}
Hence, from the bounded convergence theorem it holds that
\begin{align*}
|\nabla  p_{t}^a(x)|\lesssim|x|\int^\infty_0 r^{-(d+2)/2}\e^{-\lambda |x|^2/r}\mu_{S_t}(\dif r),
\end{align*}
and the integral  expression in the r.h.s. precisely corresponds to the stable heat kernel in dimension $d+2$ at time $t$ and point $\tilde x\in \R^{d+2} $ s.t. $|\tilde x|=\sqrt{\lambda}|x| $. Thus, from \eqref{Es0},
\begin{align*}
|\nabla  p_{t}^a(x)|&\lesssim |x| t^{-(d+2)/\alpha}\frac{1}{(1+t^{-1/\alpha} |\tilde x|)^{d+2+\alpha}}\\
&\lesssim t(t^{1/\alpha}+|x|)^{-d-\alpha-1}=\varrho^{(\alpha+1)}_{0,\alpha}(t,x).
\end{align*}
The approach is similar for higher order derivatives. This is also the case for \eqref{Es22} recalling that differentiating a Gaussian density w.r.t. the variance does not induce additional singularities.
The proof is complete.
\end{proof}
\br\rm \label{IMP_REM}
We would like to emphasize that the gradient  estimate \eqref{Es2} plays a crucial role for two-sided estimates due to the fact
that for any $\beta\in[0,1]$,
$$
|x|^\beta\varrho^{(\alpha+1)}_{0,\alpha}(t,x)=\frac{t|x|^\beta}{(t^{1/\alpha}+|x|)^{d+\alpha+1}}
\leq\frac{t^{(\alpha+\beta-1)/\alpha}}{(t^{1/\alpha}+|x|)^{d+\alpha}}=\varrho^{(\alpha)}_{0,\beta+\alpha-1}(t,x).
$$
In particular, for any $\beta\in[0,1]$,
\begin{align}\label{AA1}
|x-\theta_{t,s}(x)|^\beta\phi^{(\alpha+1)}_{0,\alpha}(s,x,t,y)
\leq\phi^{(\alpha)}_{0,\beta+\alpha-1}(s,x,t,y).
\end{align}
We carefully point out that the Gradient estimate \eqref{Es2}, which remarkably emphasizes a concentration gain, does not hold for a general $\alpha$-stable like process \cite{DZ}. This is also why, for the driving process in \eqref{SDE0}, we limit 
ourselves to the rotationally invariant, and thus symmetric,
$\alpha$-stable process and do not handle general $\alpha$-stable like processes.
\er

The following lemma is taken from \cite[Lemmas 3.2 and 3.3 ]{Ch-Zh1}.
\bl \label{Le210}
Under \eqref{AA4},  there is a constant $C=C(d,\alpha,\kappa_1)>0$ such that
\begin{align}\label{CTR_DIFF_SIGMA_GRAD}
|\nabla p_{s,t}^a-\nabla p_{s,t}^{\bar a}|(x)\lesssim_C \|a-\bar a\|_\infty\varrho^{(\alpha+1)}_{0,\alpha}(t-s, x).
\end{align}
Also,
\begin{align}\label{AA8}
|\cD^{(\alpha)}p_{s,t}^a|(x)\lesssim_C\varrho^{(\alpha)}_{0,0}(t-s, x),
\end{align}
and
\begin{align}\label{AA5}
|\cD^{(\alpha)}(p_{s,t}^a-p_{s,t}^{\bar a})|(x)\lesssim_C\|a-\bar a\|_\infty\varrho^{(\alpha)}_{0,0}(t-s, x).
\end{align}
Moreover, we also have
\begin{align}\label{AA6}
\begin{split}
&\int_{\mR^d}|\delta^{(2)}_{p_{s,t}^a}(x_1;z)-\delta^{(2)}_{p_{s,t}^a}(x_2;z)|\frac{\dif z}{|z|^{d+\alpha}}
\\&\qquad\lesssim_C\left(\frac{|x_1-x_2|}{(t-s)^{1/\alpha}}\wedge 1\right)\Bigg(\sum_{i=1,2}\varrho^{(\alpha)}_{0,0}(t-s, x_i)\Bigg).
\end{split}
\end{align}
\el
\begin{proof}
From the scaling property \eqref{SC0}, it suffices to consider $s=0$ and $t=1$. Note that
\begin{align*}
|\delta^{(2)}_{p_1^a}(x;z)|&=|p_1^a(x+z)+p_1^a(x-z)-2p_1^a(x)|\\
&\lesssim(|z|^2\wedge 1)(\varrho^{(\alpha)}(x+z)+\varrho^{(\alpha)}(x-z)+\varrho^{(\alpha)}(x)).
\end{align*}
By elementary calculations, one sees that
\begin{align}\label{AA86}
\int_{\mR^d}\varrho^{(\alpha)}(x+z)\frac{(|z|^2\wedge 1)\dif z}{|z|^{d+\alpha}}\lesssim_C\varrho^{(\alpha)}(x).
\end{align}
Thus \eqref{AA8} follows. As for \eqref{AA5} and \eqref{AA6}, they can be derived similarly 
to \cite[Lemma 2.7 and Lemma 2.8]{Ch-Zh20}. The statement \eqref{CTR_DIFF_SIGMA_GRAD} can also be derived
from the arguments developed therein. 
We omit the details.
\end{proof}

\section{Heat kernel of nonlocal operators with smooth coefficients}
\label{SEC_HK_SMOOTH_COEFF}
In this section we assume that {\bf (H$^\beta_b$)} and {\bf (H$^{\gamma}_a$)} hold and additionally that
\begin{align}\label{DA1}
\blue\|\nabla^jb\|_\infty+\|\nabla^ja\|_\infty<\infty,\ \ j\in\{0\}\cup\mN.
\end{align}
{\blue We again emphasize that we here assume that the coefficients are smooth and the drift is bounded. This last point precisely allows to derive that Duhamel like expansions holds for the semigroup (see equations \eqref{MILD_DENS_SMOOTH} and \eqref{MILD_DENS_SMOOTH0} below). We will then first remove the smoothness assumption in Subsection \ref{SUB_SEC_BD} and the boundedness assumption on the drift in Subsection \ref{SUB_SEC_TRUNC}}.

We shall denote 
$$
\sC:=\big\{(b,a): \mbox{ satisfying {\bf (H$^\beta_b$)}, {\bf (H$^{\gamma}_a$)} with common bounds $\kappa_0,\kappa_1$ and \eqref{DA1}.}\big\}
$$
Under {\bf (H$^{\gamma}_a$)} and \eqref{DA1}, 
for each $(s,x)\in\mR_+\times\mR^d$,
it is well known that there is a unique solution $(X_{s,t}(x))_{t\geq s}$ to SDE \eqref{SDE0},
and $X_{s,t}(x)$ has for $t>s$ a density $p(s,x,t,y)$ so that (cf. \cite{De-Fo, CHZ})
$$
P_{s,t} f(x):=\mE f(X_{s,t}(x))=\int_{\mR^d}f(y)p(s,x,t,y)\dif y,\ \ f\in L^\infty(\mR^d).
$$
The density is also a mild solution of the Kolmogorov equation in the sense that for all $\varphi\in C_0^2(\R^d) $
\begin{equation}
\label{MILD_DENS_SMOOTH}
P_{s,t} \varphi(x)=\varphi(x)+\int_s^t P_{s,r}\mathscr L_r \varphi (x)\dif r.
\end{equation}
{\blue Moreover, by Schauder's estimate for nonlocal parabolic equations \cite[Theorem 3.5]{Zh-Zh}, we also have $P_{s,t}\varphi\in C^\infty_b(\mR^d)$, and
$P_{s,t}\varphi$ solves the following backward Kolmogorov equation
\begin{equation}
\label{MILD_DENS_SMOOTH0}
P_{s,t} \varphi(x)=\varphi(x)+\int_s^t\mathscr L_r P_{r,t}\varphi (x)\dif r.
\end{equation}
}

Fix $(\tau,\xi)\in[s,t]\times\mR^d$. Consider the following freezing process
$$
X^{(\tau,\xi)}_{s,t}:=x+\int^t_s b_{|r-\tau|^{1/\alpha}}(r, \theta_{\tau,r}(\xi))\dif r+\int^t_s a(r, \theta_{\tau,r}(\xi))\dif L^{(\alpha)}_r.
$$
By Lemma \ref{Le27}, the density of $X^{(\tau,\xi)}_{s,t}$ is given by
\begin{align}\label{KP8}
\tilde p^{(\tau,\xi)}(s,x,t,y)= p_{s,t}^{a^{(\tau,\xi)}}\left(x-y+\int^t_s b_{|r-\tau|^{1/\alpha}}(r, \theta_{\tau,r}(\xi))\dif r\right),
\end{align}
where $a^{(\tau,\xi)}(r):=a(r,\theta_{\tau,r}(\xi)) $ and $p_{s,t}^{a^{(\tau,\xi)}}$ is the density of $\int^t_s a^{(\tau,\xi)}(r)\dif L^{(\alpha)}_r$
given in Lemma \ref{Le27}. In particular,
\begin{equation}\label{KOLMO_FREEZE}
\p_s \tilde p^{(\tau,\xi)}(s,x,t,y)+{\tilde \sL}^{(\tau,\xi)}_s{\tilde p}^{(\tau,\xi)}(s,\cdot,t,y)(x)=0,
\end{equation}
where
$$
\tilde \sL^{(\tau,\xi)}_{s} f(x):=\tfrac 12\tilde \cL^{(\tau,\xi)}_{s} f(x)+b_{|s-\tau|^{1/\alpha}}(s, \theta_{\tau,s}(\xi))\cdot\nabla f(x)
$$
and
$$
\tilde \cL^{(\tau,\xi)}_{s} f(x)=\int_{\mR^d}\delta^{(2)}_f(x;z)\frac{\kappa(s,\theta_{\tau,s}(\xi),z)}{|z|^{d+\alpha}}\dif z
$$
with
$$
\kappa(s,\theta_{s,\tau}(\xi),z):=\frac{\det(a^{-1}(s,\theta_{\tau,s}(\xi))|z|^{d+\alpha}}{|a^{-1}(s,\theta_{\tau,s}(\xi))z|^{d+\alpha}}.
$$
For simplicity, we shall write
\begin{align}
\sA^{(\tau,\xi)}_{s}f(x):=(\sL_{s}-\tilde \sL^{(\tau,\xi)}_{s})f(x)=\sK^{(\tau,\xi)}_{s}f(x)+\sB^{(\tau,\xi)}_{s}f(x),\label{DEF_A}
\end{align}
where\label{WHERE_K_B_DEFINED}
\begin{align}\label{SM1}
\sK^{(\tau,\xi)}_{s}f(x):=\tfrac 12(\cL_{s}-\tilde \cL^{(\tau,\xi)}_{s})f(x),
\end{align}
and
\begin{align}\label{SM2}
\sB^{(\tau,\xi)}_{s}f(x):=\big(b(s,x)-b_{|s-\tau|^{1/\alpha}}(s, \theta_{\tau,s}(\xi))\big)\cdot\nabla f(x).
\end{align}
Let us introduce the corresponding frozen semi-group:
\begin{align}\label{TPT}
\wt P^{(\tau,\xi)}_{s,t}f(x):=\mE f(X^{(\tau,\xi)}_{s,t}(x)).
\end{align}

We have the following Duhamel type representation formula:
\bl\label{Le310}
For any $f\in C_b^\infty(\mR^d)$ and $(\tau,\xi)\in[s,t]\times\mR^d$, it holds that
$$
P_{s,t}f=\widetilde P^{(\tau,\xi)}_{s,t}f+\int^t_sP_{s,r}\sA^{(\tau,\xi)}_r\widetilde P^{(\tau,\xi)}_{r,t}f\dif r
=\widetilde P^{(\tau,\xi)}_{s,t}f+\int^t_s\widetilde P^{(\tau,\xi)}_{s,r}\sA^{(\tau,\xi)}_rP_{r,t}f\dif r.
$$
\el
\begin{proof}
We drop for the proof the superscript $(\tau,\xi)$ for notational simplicity. From \eqref{MILD_DENS_SMOOTH} and \eqref{KOLMO_FREEZE},
$$
\p_t P_{s,t}f=P_{s,t}\sL_t f,\ \ \p_s \widetilde P_{s,t}f=-\widetilde \sL_s\widetilde P_{s,t}f.
$$
{\green Note that \blue $f\in C_b^\infty$, $\partial_r \widetilde P_{r,t}f\in C^\infty_b$}. \textcolor{black}{This indeed follows from the smoothness of $f$ and equation \eqref{KOLMO_FREEZE} through integration by parts and using the self-adjoint property of the frozen non-local operator}. We thus have by the chain rule, 
$$
\p_r(P_{s,r}\widetilde P_{r,t}f)=P_{s,r}\sL_r\widetilde P_{r,t}f-P_{s,r}\wt\sL_r\widetilde P_{r,t}f
=P_{s,r}\sA_r\widetilde P_{r,t}f.
$$
Integrating both sides from $s$ to $t$ with respect to $r$ yields
$$
P_{s,t}f=\widetilde P_{s,t}f+\int^t_sP_{s,r}\sA_r\widetilde P_{r,t}f\dif r.
$$
Similarly, by \eqref{MILD_DENS_SMOOTH0}, one can show that 
$$
\wt P_{s,t}f= P_{s,t}f-\int^t_s\widetilde P_{s,r}\sA_rP_{r,t}f\dif r.
$$
The proof is complete.
\end{proof}

By Lemma \ref{Le310}, we have for each $(\tau,\xi)\in[s,t]\times\mR^d$ and $x,y\in\mR^d$,
$$
p(s,x,t,y)
=\wt p^{(\tau,\xi)}(s,x,t,y)+
\int^t_s\int_{\R^d}p(s,x,r,z)\sA^{(\tau,\xi)}_{r}\wt p^{(\tau,\xi)}(r,\cdot,t,y)(z)\dif z\dif r.
$$
In particular, if we take $(\tau,\xi)=(t,y)$ and define
\begin{align}\label{Fre}
p_0(s,x,t,y):=\tilde p^{(t,y)}(s,x,t,y)= p_{{s,t}}^{a^{(t,y)}}\left(x-\theta_{t,s}(y)\right),
\end{align}
then  we obtain the {\it forward} representation, 
\begin{align}\label{PARAM_FORWARD}
p(s,x,t,y)
=p_0(s,x,t,y)+
\int^t_s\int_{\R^d}p(s,x,r,z)\sA^{(t,y)}_{r}p_0(r,\cdot,t,y)(z)\dif z\dif r.
\end{align}
Let 
$$
q_0(s,x,t,y):=\sA^{(t,y)}_{s}p_0(s,\cdot,t,y)(x),
$$
and define recursively for $n\ge 1 $,
\begin{equation}\label{DEF_Q}
q_n:=q_0\otimes q_{n-1},\ \ q=\sum_{n=0}^\infty q_n.
\end{equation}
By iteration, we \textit{formally} obtain from \eqref{PARAM_FORWARD} and \eqref{DEF_Q},
\begin{align}\label{Eq1}
p=p_0+p\otimes q_0=p_0+\sum_{n=0}^\infty p_0\otimes q_n=p_0+p_0\otimes q.
\end{align}

The following lemma is a direct consequence of \eqref{Fre}, \eqref{AA8} and \eqref{Es2}.
\bl \label{LEMMA_BD_PROXY}
For any $\alpha\in(0,2)$ and $j=0,1,\cdots$, we have
\begin{align}\label{Eq19}
|\nabla^j p_0(s,\cdot,t,y)|(x)\lesssim \phi^{(\alpha+j)}_{0,\alpha}(s,x,t,y)
\end{align}
and
\begin{align}\label{Eq20}
|\cD^{(\alpha)} p_0(s,\cdot,t,y)|(x)\lesssim \phi^{(\alpha)}_{0,0}(s,x,t,y).
\end{align}
\el
The following lemma corresponds to \cite[Theorem 3.1]{Ch-Zh}.
\bl\label{Le31}
The series $q=\sum_{n=0}^\infty q_n$
is absolutely convergent, and for each $s<t$, $(x,y)\mapsto q(s,x,t,y)$
is equi-continuous in $(b,a)\in\sC$. Moreover, for any $T>0$, there is a constant $C=C(T,\Theta)>0$
such that for all $(s,x,t,y)\in\mD_T$,
\begin{align}\label{G7}
|q(s,x,t,y)|\lesssim_C\big(\phi^{(\alpha)}_{\gamma_0,0}+\phi^{(\alpha)}_{0,\gamma_0}\big)(s,x,t,y),
\end{align}
where $\gamma_0:=(\alpha+\beta-1)\wedge\gamma$, and for any $\gamma_1\in(0,\gamma_0)$, 
\begin{align}\label{G8}
\begin{split}
&|q(s,x,t,y)-q(s,x',t,y)|\lesssim_C(|x-x'|^{\gamma_1}\wedge1)
\\&\quad\times\Big(\big(\phi^{(\alpha)}_{\gamma_0,-\gamma_1}+\phi^{(\alpha)}_{0,\gamma_0-\gamma_1}\big)(s,x,t,y)
+\big(\phi^{(\alpha)}_{\gamma_0,-\gamma_1}+\phi^{(\alpha)}_{0,\gamma_0-\gamma_1}\big)(s,x',t,y)\Big).
\end{split}
\end{align}
\el
\begin{proof}
(i) First of all, note that by \eqref{AA2},
\begin{align*}
|\kappa(s,x,z)-\kappa(s, \theta_{t,s}(y),z)|\lesssim (|x-\theta_{t,s}(y)|^{\gamma}\wedge 1)
\end{align*}
and by \eqref{BB0},
\begin{align*}
|b(s,x)-b_{|s-t|^{1/\alpha}}(s, \theta_{t,s}(y))|
\lesssim |x-\theta_{t,s}(y)|^\beta+|x-\theta_{t,s}(y)|+|t-s|^{\beta/\alpha}.
\end{align*}
Thus, we have by \eqref{Eq20},
\begin{align*}
|\sK^{(t,y)}_{s}p_0(s,\cdot,t,y)(x)|\lesssim \phi^{(\alpha)}_{\gamma,0}(s,x,t,y),
\end{align*}
and by \eqref{Eq19} and \eqref{AA1},
\begin{align*}
|\sB^{(t,y)}_{s} p_0(s,\cdot,t,y)(x)|\lesssim \phi^{(\alpha)}_{0,\alpha+\beta-1}(s,x,t,y).
\end{align*}
So, for $\gamma_0=\gamma\wedge(\alpha+\beta-1)$,
\begin{align*}
|q_0(s,x,t,y)|&\lesssim\Big(\phi^{(\alpha)}_{\gamma,0}+\phi^{(\alpha)}_{0,\alpha+\beta-1}\Big)(s,x,t,y)
\lesssim\Big(\phi^{(\alpha)}_{\gamma_0,0}+\phi^{(\alpha)}_{0,\gamma_0}\Big)(s,x,t,y),
\end{align*}
\textcolor{black}{recalling \eqref{DOM_STAB} for the last inequality.}
Suppose now that for some $k\in\mN$,
$$
|q_{k-1}(s,x,t,y)|\leq C_k\Big(\phi^{(\alpha)}_{\gamma_0,(k-1)\gamma_0}+\phi^{(\alpha)}_{0,k\gamma_0}\Big)(s,x,t,y).
$$
By Lemma \ref{Le09}, we have
\begin{align}
|q_{k}(s,x,t,y)|
&\leq C C_k\big(\phi^{(\alpha)}_{\gamma_0,0}+\phi^{(\alpha)}_{0,\gamma_0}\big)\otimes
 \big(\phi^{(\alpha)}_{\gamma_0,(k-1)\gamma_0}+\phi^{(\alpha)}_{0,k\gamma_0}\big)(s,x,t,y)
\notag\\
&\leq C_0C_k\cB(\tfrac{\gamma_0}{\alpha},\tfrac{k\gamma_0}{\alpha})
 \big(\phi^{(\alpha)}_{\gamma_0,k\gamma_0}+\phi^{(\alpha)}_{0,(k+1)\gamma_0}\big)(s,x,t,y).\label{BD_Q_K}
\end{align}
Hence,
$$
C_{k+1}=C_0C_k\cB(\tfrac{\gamma_0}{\alpha},\tfrac{k\gamma_0}{\alpha}).
$$
From the relation $\cB(\gamma,\beta)=\frac{\Gamma(\gamma)\Gamma(\beta)}{\Gamma(\gamma+\beta)}$, where $\Gamma$
is the usual Gamma function, we obtain
$$
C_k=C^k_0\prod_{i=1}^{k-1}\cB(\tfrac{\gamma_0}{\alpha},\tfrac{(k-1)\gamma_0}{\alpha})
=\frac{(C_0\Gamma(\gamma_0/\alpha))^k}{\Gamma(k\gamma_0/\alpha)},
$$
with the usual convention that $\prod_{i=1}^{0}=1 $.
Thus
\begin{align*}
\sum_{k=0}^\infty |q_k(s,x,t,y)|
&\leq \sum_{k=0}^\infty\frac{(C_0\Gamma(\gamma_0/\alpha))^k}{\Gamma(k\gamma_0/\alpha)}
\Big(\phi^{(\alpha)}_{\gamma_0,k\gamma_0}+\phi^{(\alpha)}_{0,(k+1)\gamma_0}\Big)(s,x,t,y)
\\&\leq \sum_{k=0}^\infty\frac{(C_0\Gamma(\gamma_0/\alpha))^k}{\Gamma(k\gamma_0/\alpha)}
\Big(\phi^{(\alpha)}_{\gamma_0,0}+\phi^{(\alpha)}_{0,\gamma_0}\Big)(s,x,t,y).
\end{align*}
This gives \eqref{G7}.\\
(ii) For fixed $s<t$, by Lemma \ref{Le27} and the definition of $q_0$, one sees that $(x,y)\mapsto q_0(s,x,t,y)$ is equi-continuous in $(b,a)\in\sC$. 
Furthermore, it follows by induction that, for each $k\in\mN$,
$(x,y)\mapsto q_k(s,x,t,y)$ is also  equi-continuous in $(b,a)\in\sC$. Hence, 
$(x,y)\mapsto q(s,x,t,y)$ is equi-continuous in $(b,a)\in\sC$. 
\smallskip\\
(iii) If $|x-x'|\geq (t-s)^{1/\alpha}$, then we have
\begin{align*}
|q_0(s,x,t,y)|&\lesssim (|x-x'|^{\gamma_1}\wedge 1)(t-s)^{-\gamma_1/\alpha}\big(\phi^{(\alpha)}_{\gamma_0,0}+\phi^{(\alpha)}_{0,\gamma_0}\big)(s,x,t,y)
\\&=(|x-x'|^{\gamma_1}\wedge 1)\big(\phi^{(\alpha)}_{\gamma_0,-\gamma_1}+\phi^{(\alpha)}_{0,\gamma_0-\gamma_1}\big)(s,x,t,y).
\end{align*}
Next we assume $|x-x'|\leq (t-s)^{1/\alpha}$. In this case,  it is easy to see from \eqref{DEF_VAR_RHO}-\eqref{DEF_PHI}, that
\begin{align}\label{AA20}
\phi^{(\eta)}_{0,0}(s,x,t,y)\asymp\phi^{(\eta)}_{0,0}(s,x',t,y),\ \ \eta\geq 0.
\end{align}
By \eqref{AA8}, \eqref{AA6} and \eqref{AA20}, we have
\begin{align*}
&|\sK^{(t,y)}_{s}p_0(s,\cdot,t,y)(x)-\sK^{(t,y)}_{s}p_0(s,\cdot,t,y)(x')|
\\&\quad\leq \|\kappa(\cdot,x,\cdot)-\kappa(\cdot,x',\cdot)\|_\infty \int_{\mR^d}|\delta^{(2)}_{p_0(s,\cdot,t,y)}(x;z)|\frac{\dif z}{|z|^{d+\alpha}}
\\&\qquad+\|\kappa(\cdot,x,\cdot)-\kappa(\cdot,\theta_{t,s}(y),\cdot)\|_\infty
\\&\qquad\qquad\times\int_{\mR^d}|\delta^{(2)}_{p_0(s,\cdot,t,y)}(x;z)-\delta^{(2)}_{p_0(s,\cdot,t,y)}(x';z)|\frac{\dif z}{|z|^{d+\alpha}}
\\&\quad\leq (|x-x'|^{\gamma}\wedge 1)\phi^{(\alpha)}_{0,0}(s,x,t,y)+(|x-\theta_{t,s}(y)|^\gamma\wedge 1)
\\&\qquad\times\left(\frac{|x-x'|}{(t-s)^{1/\alpha}}\wedge 1\right)\Big(\phi^{(\alpha)}_{0,0}(s,x,t,y)+\phi^{(\alpha)}_{0,0}(s,x',t,y)\Big)
\\&\quad\lesssim (|x-x'|^{\gamma_1}\wedge 1)\Big(\phi^{(\alpha)}_{0,\gamma-\gamma_1}(s,x,t,y)+\phi^{(\alpha)}_{\gamma,-\gamma_1}(s,x,t,y)\Big)\textcolor{black}{.}
\end{align*}
Moreover, by \eqref{Eq19}, \eqref{AA20} and \eqref{AA1}, we also have
\begin{align*}
&|\sB^{(t,y)}_{s} p_0(s,\cdot,t,y)(x)-\sB^{(t,y)}_{s} p_0(s,\cdot,t,y)(x')|
\leq |b(s,x)-b(s,x')|\cdot|\nabla p_0(s,\cdot,t,y)|(x')
\\&\qquad+\big|b(s,x)-b_{|s-t|^{1/\alpha}}(s, \theta_{t,s}(y))\big|
\cdot |\nabla p_0(s,\cdot,t,y)(x')-\nabla p_0(s,\cdot,t,y)(x)|
\\&\quad\lesssim|x-x'|^\beta\phi^{(\alpha+1)}_{0,\alpha}(s,x,t,y)+(|x-\theta_{t,s}(y)|^\beta+|t-s|^{\beta/\alpha})|x-x'|\phi^{(\alpha+2)}_{0,\alpha}(s,x,t,y)
\\&\quad\lesssim(|x-x'|^{\gamma_1}\wedge 1)\phi^{(\alpha)}_{0,\alpha+\beta-1-\gamma_1}(s,x,t,y).
\end{align*}
Combining the above calculations and recalling $\gamma_0=\gamma\wedge(\alpha+\beta-1)$, we obtain
\begin{align*}
&|q_0(s,x,t,y)-q_0(s,x',t,y)|\lesssim_C(|x-x'|^{\gamma_1}\wedge1)
\\&\quad\times\Big(\big(\phi^{(\alpha)}_{\gamma_0,-\gamma_1}+\phi^{(\alpha)}_{0,\gamma_0-\gamma_1}\big)(s,x,t,y)
+\big(\phi^{(\alpha)}_{\gamma_0,-\gamma_1}+\phi^{(\alpha)}_{0,\gamma_0-\gamma_1}\big)(s,x',t,y)\Big).
\end{align*}
Using this last estimate, equation \eqref{G8} follows from the same iterative argument as in  (i).
\end{proof}
\br\rm
This lemma allows to iterate the representation \eqref{PARAM_FORWARD} which leads to the representation \eqref{Eq1} of the density.
\er
 We now aim at proving the following \textit{a priori} estimate about $p(s,x,t,y)$.
\bt\label{T8} 
Under {\bf (H$^\gamma_a$)}, {\bf (H$^\beta_b$)} and \eqref{DA1}, for each $0\leq s<t<\infty$,
$X_{s,t}(x)$ admits a  density $p(s,x,t,y)$ that is equi-continuous in $(b,a)\in\sC$ as a function of $x,y\in\mR^d$, 
and there is a constant $C=C(T,\Theta)>0$ so that for all $(s,x,t,y)\in\mD_T$,
\begin{align}\label{TW0}
p(s,x,t,y)\asymp_C \phi^{(\alpha)}_{0,\alpha}(s,x,t,y).
\end{align}
\et
\begin{proof}
Note that by \eqref{Fre}, \eqref{Es0} and \eqref{KH2},
$$
p_0(s,x,t,y)\asymp_C\phi^{(\alpha)}_{0,\alpha}(s,x,t,y).
$$
By Lemma \ref{Le09}, we have
\begin{align*}
|p_0\otimes q|(s,x,t,y)
\lesssim_C(\phi^{(\alpha)}_{0,\alpha+\gamma_0}+\phi^{(\alpha)}_{\gamma_0,\alpha})(s,x,t,y).
\end{align*}
The upper bound follows from \eqref{Eq1}.

Next we use Lemma \ref{Le26} to show the lower bound estimate. Let $K$ be as in  Lemma \ref{Le26}.
Suppose that $|x-\theta_{t,s}(y)|\leq 2K(t-s)^{1/\alpha}$ \textcolor{black}{(diagonal regime)}. Then we have
\begin{align*}
p(s,x,t,y)&\geq p_0(s,x,t,y)-|p_0\otimes q(s,x,t,y)|
\\&\geq c_0\phi^{(\alpha)}_{0,\alpha}(s,x,t,y)-(\phi^{(\alpha)}_{0,\alpha+\gamma_0}+\phi^{(\alpha)}_{\gamma_0,\alpha})(s,x,t,y)
\\&\geq (c_0-C_1(t-s)^{\frac{\gamma_0}{\alpha}})\phi^{(\alpha)}_{0,\alpha}(s,x,t,y),
\end{align*}
\textcolor{black}{recalling from \eqref{DEF_PHI} and \eqref{DEF_VAR_RHO} that, in the current diagonal regime, $ (\phi^{(\alpha)}_{0,\alpha+\gamma_0}+\phi^{(\alpha)}_{\gamma_0,\alpha})(s,x,t,y)\le C_1(t-s)^{\frac{\gamma_0}{\alpha}}\phi^{(\alpha)}_{0,\alpha}(s,x,t,y)$ for the last inequality}.
In particular, if $t-s\leq\ell$ with $\ell$ small enough and $|x-\theta_{t,s}(y)|\leq 2K(t-s)^{1/\alpha}$, then
\begin{align}\label{DE1}
p(s,x,t,y)\geq \tfrac{c_0}{2}\phi^{(\alpha)}_{0,\alpha}(s,x,t,y)\geq c_1(t-s)^{-d/\alpha}.
\end{align}
Next we prove the above estimate still holds for 
$$
|x-\theta_{t,s}(y)|\geq 2K(t-s)^{1/\alpha}.
$$
Let $\eps\in(0,1/2)$ be as in Lemma \ref{Le26} and small enough 
so that $2(1-\eps)^{1/\alpha}\geq 1$. Let 
$$
r:=s+\eps(t-s),\ \ B:=\{z: |z-\theta_{t,r}(y)|\leq 2K(t-r)^{1/\alpha}\}.
$$
Since $2(1-\eps)^{1/\alpha}\geq 1$, we clearly have
$$
B\supset\{z: |z-\theta_{t,r}(y)|\leq K(t-s)^{1/\alpha}\}=:B'.
$$
Now from the Chapman-Kolmogorov equation, we have for $t-s\leq\ell$,
\begin{align*}
p(s,x,t,y)&=\int_{\mR^d}p(s,x,r,z)p(r,z,t,y)\dif z
\\&\geq\int_B p(s,x,r,z)p(r,z,t,y)\dif z
\\&\geq\inf_{z\in B}p(r,z,t,y)\int_B p(s,x,r,z)\dif z
\\&\!\!\!\stackrel{\eqref{DE1}}{\geq} c_1(t-r)^{-d/\alpha}\mP(X_{s,r}(x)\in B)
\\&\geq c_2(t-s)^{-d/\alpha}\mP(X_{s,r}(x)\in B')
\\&\geq c_3(t-s)|x-\theta_{t,s}(y)|^{-d-\alpha},
\end{align*}
where the last step is due to Lemma \ref{Le26}.
Thus we obtain that for some $c_4>0$ and all $s,t\in[0,T]$,
$$
p(s,x,t,y)\geq c_4\phi^{(\alpha)}_{0,\alpha}(s,x,t,y),\  t-s\leq\ell, \ x,y\in\mR^d.
$$
For $t-s>\ell$, the bound follows iteratively from the Chapman-Kolmogorov equation.
The proof is complete.
\end{proof}
For the fractional derivative estimates, we need the following lemma.
\bl\label{Le34}
For $s<t$, let $h_{s,t}(x):=\int_{\mR^d}p_0(s,x,t,y)\dif y$. We have for some $C>0$,
$$
|\cD^{(\alpha)}h_{s,t}|(x)\lesssim_C(t-s)^{\gamma_0/\alpha-1},\  \gamma_0:=\gamma \wedge (\alpha+\beta-1).
$$
\el
\begin{proof}
By definition we have
\begin{align*}
|\cD^{(\alpha)}h_{s,t}|(x)
&=\int_{\mR^d}\left|\int_{\mR^d}\delta^{(2)}_{\tilde p^{(t,y)}(s,\cdot,t,y)}(x;z)\dif y\right|\frac{\dif z}{|z|^{d+\alpha}}
\\&=\int_{\mR^d}\left|\int_{\mR^d}\delta^{(2)}_{p_{{s,t}}^{a^{(t,y)}}}(x-\theta_{t,s}(y);z)\dif y\right|\frac{\dif z}{|z|^{d+\alpha}}
\\&\leq\int_{\mR^d}\left|\int_{\mR^d}\delta^{(2)}_{p_{{s,t}}^{a^{(t,y)}}-p_{{s,t}}^{a^{(s,x)}}} (x-\theta_{t,s}(y);z)\dif y\right|\frac{\dif z}{|z|^{d+\alpha}}
\\&+\int_{\mR^d}\left|\int_{\mR^d}\delta^{(2)}_{p_{{s,t}}^{a^{(s,x)}}}(x-\theta_{t,s}(y);z)\dif y\right|\frac{\dif z}{|z|^{d+\alpha}}=:I_1+I_2.
\end{align*}
For $I_1$, noting that by \eqref{AA0} and Lemma \ref{Le22},
$$
|a(r,\theta_{s,r}(x))-a(r,\theta_{t,r}(y))|\lesssim 1\wedge |x-\theta_{t,s}(y)|^{\gamma}+|t-s|^{\gamma/\alpha},
$$
we have
\begin{align*}
I_1&\leq\int_{\mR^d}|\cD^{(\alpha)}(p_{s,t}^{a^{(t,y)}}-p_{s,t}^{a^{(s,x)}})|(x-\theta_{t,s}(y))\dif y
\\&\stackrel{\eqref{AA5}}{\lesssim}\int_{\mR^d}\Big(\phi^{(\alpha)}_{\gamma,0}+\phi^{(\alpha)}_{0,\gamma}\Big)(s,x,t,y)\dif y
\stackrel{\eqref{ES4}}{\lesssim} (t-s)^{\gamma/\alpha-1}.
\end{align*} 
For $I_2$, by the change of variable we have
\begin{align*}
I_2&=\int_{\mR^d}\left|\int_{\mR^d}\delta^{(2)}_{p_{s,t}^{a^{(s,x)}}}(x-y;z)\det(\nabla\theta^{-1}_{s,t}(y))\dif y\right|\frac{\dif z}{|z|^{d+\alpha}}
\\&=\int_{\mR^d}\left|\int_{\mR^d}\delta^{(2)}_{p_{s,t}^{a^{(s,x)}}}(x-y;z)\Big(\det(\nabla\theta^{-1}_{s,t}(y))-1\Big)\dif y\right|\frac{\dif z}{|z|^{d+\alpha}},
\end{align*}
where we have used that
$$
\int_{\mR^d}p_{s,t}^{a^{(s,x)}} (x-y)\dif y=1\Rightarrow \int_{\mR^d}\delta^{(2)}_{p_{s,t}^{a^{(s,x)}}}(x-y)\dif y=0.
$$
Thus by (i) of Lemma \ref{Le22} and \eqref{AA8}, we have
\begin{align*}
I_2&\lesssim(t-s)^{\frac{\beta+\alpha-1}{\alpha}}\int_{\mR^d}|\cD^{(\alpha)}p_{s,t}^{a^{(s,x)}}|(x-y)\dif y
\\&\lesssim(t-s)^{\frac{\beta+\alpha-1}{\alpha}}\int_{\mR^d}\varrho^{(\alpha)}_{0,0}(t-s,x-y)\dif y\\
&\lesssim(t-s)^{\frac{\beta-1}{\alpha}}=(t-s)^{-1+\frac{\alpha+\beta-1}{\alpha}}.
\end{align*}
The proof is complete.
\end{proof}
\bl\label{Le35}
(Fractional derivative estimate) For any $T>0$,  we have for some $C=C(T,\Theta)>0$,
$$
|\cD^{(\alpha)}p(s,\cdot,t,y)|(x)\lesssim_C \phi^{(\alpha)}_{0,0}(s,x,t,y).
$$
\el
\begin{proof}
Let $u=(s+t)/2$. By \eqref{Eq1} and the definition of $\delta^{(2)}$, we have
\begin{align*}
\delta^{(2)}_{p(s,\cdot,t,y)}(x;\bar z)&=
\delta^{(2)}_{p_0(s,\cdot,t,y)}(x;\bar z)+\int^t_{s}\!\!\!\int_{\mR^d}\delta^{(2)}_{p_0(s,\cdot,r,z)}(x;\bar z)q(r,z,t,y)\dif z\dif r\\
&=\delta^{(2)}_{p_0(s,\cdot,t,y)}(x;\bar z)+\int^u_{s}\left(\int_{\mR^d}\delta^{(2)}_{p_0(s,\cdot,r,z)}(x;\bar z)\dif z\right)q(r,\theta_{s,r}(x),t,y)\dif r\\
&+\int^u_{s}\!\!\!\int_{\mR^d}\delta^{(2)}_{p_0(s,\cdot,r,z)}(x;\bar z)(q(r,z,t,y)-q(r,\theta_{s,r}(x),t,y))\dif z\dif r\\
&+\int^t_u\!\!\!\int_{\mR^d}\delta^{(2)}_{p_0(s,\cdot,r,z)}(x;\bar z)q(r,z,t,y)\dif z\dif r.
\end{align*}
{ With the notations of Lemma \ref{Le34}, set} 
$h_{s,r}(x)=\int_{\mR^d}p_0(s,x,r,z)\dif z$. By \eqref{DA0} and the Fubini theorem, we have
\begin{align*}
|\cD^{(\alpha)} p(s,\cdot,t,y)|(x)
&\leq|\cD^{(\alpha)} p_0(s,\cdot,t,y)|(x)+\int^u_s |\cD^{(\alpha)}h_{s,r}|(x)|q(r,\theta_{s,r}(x),t,y)|\dif r
\\&+\int^u_s\!\!\!\int_{\mR^d}|\cD^{(\alpha) }p_0(s,\cdot,r,z)|(x)|q(r,z,t,y)-q(r,\theta_{s,r}(x),t,y)|\dif z\dif r
\\&+\int^t_u\!\!\!\int_{\mR^d}|\cD^{(\alpha) }p_0(s,\cdot,r,z)|(x)|q(r,z,t,y)|\dif z\dif r\\
&=:I_1(x)+I_2(x)+I_3(x)+I_4(x).
\end{align*}
For $I_1$, by \eqref{Eq20} we have
$$
I_1(x)\lesssim \phi^{(\alpha)}_{0,0}(s,x,t,y).
$$
Recall
$$
\gamma_0\textcolor{black}{=}(\alpha+\beta-1)\wedge\gamma,\quad \gamma_1\in(0,\gamma_0).
$$
For $I_2$, by Lemma \ref{Le34}, \eqref{G7}, \eqref{ES44} and \eqref{KH1}, we have
\begin{align*}
I_2(x)&\lesssim\int^u_s(r-s)^{\frac{\gamma_0}{\alpha}-1}\big(\phi^{(\alpha)}_{\gamma_0,0}
+\phi^{(\alpha)}_{0,\gamma_0}\big)(r,\theta_{s,r}(x),t,y)\dif r
\\&\lesssim\left(\int^u_s(r-s)^{\frac{\gamma_0}{\alpha}-1}\dif r\right)
\big(\phi^{(\alpha)}_{\gamma_0,0}+\phi^{(\alpha)}_{0,\gamma_0}\big)(s,x,t,y)
\\&\lesssim\big(\phi^{(\alpha)}_{\gamma_0,\gamma_0}+\phi^{(\alpha)}_{0,2\gamma_0}\big)(s,x,t,y).
\end{align*}
For $I_3$, by \textcolor{black}{\eqref{Eq20}}, \eqref{G8} and \eqref{eq31}, we have
\begin{align*}
I_3(x)&\lesssim\int^u_s\int_{\mR^d}\phi^{(\alpha)}_{\gamma_1,0}(s,x,r,z)\big(\phi^{(\alpha)}_{\gamma_0,-\gamma_1}
+\phi^{(\alpha)}_{0,\gamma_0-\gamma_1}\big)(r,z,t,y)\dif z\dif r\\
&+\int^u_s\int_{\mR^d}\phi^{(\alpha)}_{\gamma_1,0}(s,x,r,z)\dif z
\big(\phi^{(\alpha)}_{\gamma_0,-\gamma_1}
+\phi^{(\alpha)}_{0,\gamma_0-\gamma_1}\big)(s,x,t,y)\dif r\\
&\lesssim \big(\phi^{(\alpha)}_{\gamma_0,0}+\phi^{(\alpha)}_{0,\gamma_0}\big)(s,x,t,y).
\end{align*}
For $I_4$, by \eqref{Eq20}, \eqref{G7}   and \eqref{eq30}, we have
\begin{align*}
I_4(x)&\lesssim\int^t_u\big(\phi^{(\alpha)}_{0,0}\odot(\phi^{(\alpha)}_{\gamma_0,0}
+\phi^{(\alpha)}_{0,\gamma_0})\big)_r(s,x,t,y)\dif r\lesssim\phi^{(\alpha)}_{0,0}(s,x,t,y).
\end{align*}
Combining the above estimates, we complete the proof.
\end{proof}

\section{A priori gradient estimates}\label{SEC_GRAD}

The aim of this section is to show the following a priori gradient estimate.
\bt\label{Th41}
Under {\bf (H$^\beta_b$)}, {\bf (H$^{\gamma}_a$)} and \eqref{DA1}, for any $T>0$, there is a constant $C=C(T,\Theta)>0$ such that
for all $f\in \cB_b(\mR^d)$,  
$0\leq s<t\leq T$ and $x\in\mR^d$,
\begin{align}\label{DQ1}
|\nabla P_{s,t}f(x)|\lesssim_C (t-s)^{-1/\alpha}P_{s,t}|f|(x).
\end{align}
Moreover, $x\mapsto \nabla P_{s,t}f(x)$ is equi-continuous in $(b,a)\in\sC$.
\et
\textcolor{black}{We again emphasize that Theorem \ref{Th41} gives that the constants in the gradient estimates actually \textbf{do not} depend on the smoothness of the coefficients, neither on the boundedness of the drift assumed in \eqref{DA1}}.

We shall prove this theorem for $\alpha\in[1,2)$ and $\alpha\in(0,1)$ separately by different methods.

\subsection{Critical and Subcritical cases: $\alpha\in[1,2)$}
In this subsection we start from the series expansion \eqref{Eq1} for the density to derive the estimate 
\begin{align}\label{GR}
|\nabla_x p(s,x,t,y)|\lesssim_{C_3}\phi^{(\alpha)}_{0,\alpha-1}(s,x,t,y),
\end{align}
 when  {\bf (H$^\beta_b$)}, {\bf (H$^{\gamma}_a$)} and \eqref{DA1} are in force and $\alpha\in[1,2)$. 
 {This precisely gives \eqref{DQ1}}. We recall that, with the notations of Section \ref{SEC_HK_SMOOTH_COEFF}:
$$
p(s,x,t,y)=p_0(s,x,t,y)+(p_0\otimes q)(s,x,t,y).
$$
\textcolor{black}{We will first assume that $\alpha\in (1,2) $ and handle the critical case $\alpha=1 $ afterwards through a domination argument. 
Now,  for $\alpha\in (1,2) $}, $u=\frac{s+t}{2}$ and $\xi=\theta_{s,r}(x)$,
\begin{align}
\nabla_x p(s,x,t,y)&=\nabla_xp_0(s,x,t,y)
+\int_{u}^{t} \int_{\R^d}\nabla_x {p}_0(s,x,r,z)q(r,z,t,y)\dif z \dif r\notag\\
&\quad+\int_s^u\int_{\R^d}(\nabla_x {p}_0-\nabla_x \tilde{p}^{(r,\xi)})(s,x,r,z)q(r,z,t,y)\dif z \dif r  \notag\\
&\quad +  \int_s^u\int_{\R^d}\nabla_x\tilde{p}^{(r,\xi)}(s,x,r,z)(q(r,z,t,y)-q(r,\xi,t,y))\dif z \dif r\notag\\
&=:G_1(s,x,t,y)+G_2(s,x,t,y)+G_3(s,x,t,y)+G_4(s,x,t,y),\label{CUT_GRAD_SUB_CRITICAL}
\end{align}
where  for the last term, we have used precisely the cancellation property
$$
\int_{\R^d} \nabla_x\tilde{p}^{(\textcolor{black}{r},\xi)}(s,x,r,z)\dif z=0.
$$ 
For $G_1$, by \eqref{Eq19} we clearly have
$$
|G_1(s,x,t,y)|\lesssim \phi^{(\alpha+1)}_{0,\alpha}(s,x,t,y)\leq\phi^{(\alpha)}_{0,\alpha-1}(s,x,t,y),
$$
using Remark \ref{IMP_REM}, equation \eqref{AA1}, for the last inequality.
For $G_2$, by \eqref{Eq19}, \eqref{G7} and \eqref{eq31},  we have
\begin{align*}
|G_2(s,x,t,y)|
&\le \int_u^t\int_{\R^d } \phi_{0,\alpha}^{(\alpha+1)}(s,x,r,z)   |q(r,z,t,y)|\dif z\dif r\notag\\
&\lesssim (t-s)^{- \frac 1\alpha}\phi_{0,\alpha}^{(\alpha)}\otimes\big(\phi_{\gamma_0,0}^{(\alpha)} +\phi_{0,\gamma_0}^{(\alpha)}\big)(s,x,t,y)\\
&\lesssim (t-s)^{-\frac 1\alpha}  \phi_{0,\alpha+\gamma_0}^{(\alpha)}(s,x,t,y)=\phi_{0,\alpha+\gamma_0-1}^{(\alpha)}(s,x,t,y).
\end{align*}
For $G_3$, noting that by \eqref{KP8},
\begin{align*}
\nabla_x {p}_0(s,x,r,z)&=\nabla_x {p}^{a^{(r,z)}}_{s,r}\left(x-z+\int^r_s b_{|r'-r|^{1/\alpha}}(r',\theta_{r,r'}(z))\dif r'\right),\\
\nabla_x \tilde{p}^{(\textcolor{black}{r},\xi)}(s,x,r,z)&=\nabla_x {p}^{a^{(r,\xi)}}_{s,r}\left(x-z+\int^r_s b_{|r'-r|^{1/\alpha}}(r',\theta_{r,r'}(\xi))\dif r'\right),
\end{align*}
by \eqref{CTR_DIFF_SIGMA_GRAD}, \eqref{Es2}, \eqref{LL01} and \eqref{AQ3}, one finds that
\begin{align}\label{DZ9}
 |\nabla_x {p}_0-\nabla_x \tilde{p}^{(r,\xi)}|(s,x,r,z)
 &\lesssim\phi_{0,\alpha}^{(\alpha+1)}(s,x,r,z) (1\wedge |z-\theta_{s,r}(x)|^\gamma)\no\\
 &+\phi_{0,\alpha}^{(\alpha+2)}(s,x,r,z)\big(|z-\theta_{s,r}(x)|^\beta+(r-s)^{\frac{\beta}{\alpha}}\big)(r-s)\no\\
 &\underset{\eqref{AA1}}{\lesssim}\big(\phi_{0,\alpha+\gamma-1}^{(\alpha)}+\phi_{0,2\alpha+\beta-2}^{(\alpha)}\big)(s,x,r,z)\no\\
 &\lesssim\phi_{0,\alpha-1+\gamma_0}^{(\alpha)}(s,x,r,z),
\end{align}
where $\gamma_0=\gamma\wedge(\alpha+\beta-1)$.
Therefore, due to $\alpha\in[1,2)$, by \eqref{eq31},
\begin{align}\label{CTR_SECOND_SUMMAND} 
|G_3(s,x,t,y)|&\lesssim 
\phi_{0,\alpha-1+\gamma_0}^{(\alpha)}\otimes
\big(\phi_{\gamma_0,0}^{(\alpha)} +\phi_{0,\gamma_0}^{(\alpha)}\big)(s,x,t,y)\no\\
&\lesssim\phi_{0,\alpha-1+2\gamma_0}^{(\alpha)}(s,x,t,y)
 \leq\phi_{0,\alpha-1}^{(\alpha)}(s,x,t,y).
\end{align}
For $G_4$, by \eqref{KP8}, \eqref{Es2}
and \eqref{G8} we have for $\gamma_1\in(0,\gamma_0)$,
\begin{align*}
|G_4(s,x,t,y)|
&\le \int_s^u\int_{\R^d}|\nabla_x\tilde{p}^{(r,\xi)}(s,x,r,z)|\ |(q(r,z,t,y)-q(r,\xi,t,y))|\dif z \dif r\\
&\lesssim  \int_s^u\dif r\int_{\R^d} \dif z  \phi_{0,\alpha-1}^{(\alpha)}(s,x,r,z)(1\wedge |z-\xi|^{\gamma_1})\\
&\times\Big[\big(\phi^{(\alpha)}_{\gamma_0,-\gamma_1}+\phi^{(\alpha)}_{0,\gamma_0-\gamma_1}\big)(r,z,t,y)
+\big(\phi^{(\alpha)}_{\gamma_0,-\gamma_1}+\phi^{(\alpha)}_{0,\gamma_0-\gamma_1}\big)(r,\xi,t,y)\Big].
\end{align*}
Since $t-r \asymp t-s$ for $r\in[s,u]$ and $\xi=\theta_{s,r}(x)$, from \eqref{KH2} in Lemma \ref{Le22}, it holds
$$
\big(\phi^{(\alpha)}_{\gamma_0,-\gamma_1}+\phi^{(\alpha)}_{0,\gamma_0-\gamma_1}\big)(r,\xi,t,y)\lesssim \big(\phi^{(\alpha)}_{\gamma_0,-\gamma_1}+\phi^{(\alpha)}_{0,\gamma_0-\gamma_1}\big)(s,x,t,y).
$$
Thus by \eqref{eq31}, we eventually have
\begin{align*}
|G_4(s,x,t,y)|
&\lesssim \int_s^u\dif r \int_{\R^d} \dif z  \phi_{0,\alpha+\gamma_1-1}^{(\alpha)}(s,x,r,z)\notag\\
&\times\Big[ \big(\phi^{(\alpha)}_{\gamma_0,-\gamma_1}+\phi^{(\alpha)}_{0,\gamma_0-\gamma_1}\big)(r,z,t,y)+ \big(\phi^{(\alpha)}_{\gamma_0,-\gamma_1}+\phi^{(\alpha)}_{0,\gamma_0-\gamma_1}\big)(s,x,t,y) \Big]\\
&\lesssim \big(\phi^{(\alpha)}_{0,\gamma_0-\gamma_1}+\phi^{(\alpha)}_{\gamma_0,\alpha-1}\big)(s,x,t,y).
\end{align*}
Combining the above calculations, we obtain \eqref{GR} in the case $\alpha\in (1,2) $.\\
\textcolor{black}{We mention that  in the critical case $\alpha=1 $, the previous terms $G_1$ and $G_2$ in \eqref{CUT_GRAD_SUB_CRITICAL} are well-defined (no time singularity) and controlled similarly. However, some care is needed to justify that
\begin{equation}\label{LIMIT_CRITICAL}
 \nabla_x \int_s^u\int_{\R^d} {p}_0(s,x,r,z)q(r,z,t,y)\dif z \dif r =\int_s^u\int_{\R^d}\nabla_x {p}_0(s,x,r,z)q(r,z,t,y)\dif z \dif r.
 \end{equation}
The previous controls on $G_3,G_4$ can be used to  prove that setting
$$\int_{s+\varepsilon}^u\int_{\R^d}\nabla_x {p}_0(s,x,r,z)q(r,z,t,y)\dif z \dif r=:\int_{s+\eps}^u 
\Gamma_{s,x,t,y}(r)\dif r,$$ the function $r\mapsto \Gamma_{s,x,t,y}(r) $ is integrable.
Hence, a domination argument yields \eqref{LIMIT_CRITICAL} which in turns establishes that \eqref{CUT_GRAD_SUB_CRITICAL} still holds, as well as the associated estimates, in the critical case $\alpha=1 $}.

\textcolor{black}{Let us now turn to the equicontinuity part of the theorem. From}
 the dominated convergence theorem \textcolor{black}{and} the above calculations,
it is easy to see that
$$
\lim_{x\to x_0}\sup_{(b,a)\in\sC}|G_i^{b,a}(s,x,t,y)-G_i^{b,a}(s,x_0,t,y)|=0,\ { i=1,2,3,4},
$$
where $G_i^{b,a}$ are defined as above through the coefficients $b,a$. For instance,
\begin{align*}
&\lim_{x\to x_0}\sup_{(b,a)\in\sC}|G_2^{b,a}(s,x,t,y)-G_2^{b,a}(s,x_0,t,y)|\\
&\leq\int_{u}^{t}\lim_{x\to x_0} \sup_{(b,a)\in\sC}\int_{\R^d}|\nabla_x {p}_0(s,x,r,z)-\nabla_x {p}_0(s,x_0,r,z)|\,|q(r,z,t,y)|\dif z \dif r,
\end{align*}
and for each $r\in(u,t)$, by \eqref{G7} and \eqref{Eq19},
\begin{align*}
&\lim_{x\to x_0} \sup_{(b,a)\in\sC}\int_{\R^d}|\nabla_x {p}_0(s,x,r,z)-\nabla_x {p}_0(s,x_0,r,z)|\,|q(r,z,t,y)|\dif z\\
&\lesssim\int_{\R^d}\lim_{x\to x_0} \sup_{(b,a)\in\sC}|\nabla_x {p}_0(s,x,r,z)-\nabla_x {p}_0(s,x_0,r,z)|\phi^{(\alpha)}_{0,0}(r,z,t,y)\dif z=0.
\end{align*}
In particular, Theorem \ref{Th41} holds for $\alpha\in[1,2)$.
\br\rm
We remark that for $\alpha\in(0,1)$, under $\alpha+\beta>1$, the second inequality in \eqref{CTR_SECOND_SUMMAND} 
may not hold since $\alpha+\gamma_0-1$ may be less than zero.
This is also the reason that we have to make a different treatment for supercritical case.
Let us mention that this proof anyhow works even in the super-critical case under the most stringent condition $\alpha+\frac \beta 2>1 $. Eventually, we also point out that the previous arguments can be simplified if $\alpha\in (1,2)$ for which the full parametrix expansion \eqref{Eq1} of the density can actually be directly differentiated since the induced singularity in time remains integrable.
\er
\subsection{Supercritical case $\alpha\in(0,1)$}
The following gradient estimate comes in \cite{Wa-Xu-Zh}.
\bt\label{Le36}
(Gradient estimate) Under {\bf (H$^\beta_b$)}, {\bf (H$^{\gamma}_a$)} and \eqref{DA1}, 
for any $T>0$,
there is a constant $C>0$ such that
for all $f\in C_b(\mR^d)$ 
and $0\leq s<t\leq T$,
$$
|\nabla P_{s,t}f(x)|\lesssim_C (t-s)^{-1/\alpha}\|f\|_\infty,
$$
{ where the constant $C$ may depend on $\|\nabla b\|_\infty$ and $\|\nabla a\|_\infty$}.
\et
\textcolor{black}{Theorem \ref{Le36} is important since it will precisely allow to justify that the Gronwall-Volterra lemma applies in the procedure below (see (Step 5) p. \pageref{USE_THM_4.3}). It is used as a \textit{prior} estimate. The analysis we now develop actually aims at proving that the constants in the gradient estimate are indeed \textbf{independent} of the smoothness of the coefficients \textcolor{black}{and the boundedness of the drift}. This is the main achievement of Theorem \ref{Th41} which we here prove for the supercritical case $\alpha\in (0,1) $}.

Below we fix $s<t$ and $x\in\mR^d$ and {\blue assume $f\in C^\infty_b(\mR^d)$.} We divide the proof into six steps.  

{\it (Step 1).} For notational simplicity, we shall write for $r\in[s,t]$,
$$
\wt\sA_r:=\sA^{(s,x)}_r=\sK^{(s,x)}_r+\sB^{(s,x)}_r=:\wt\sK_r+\wt\sB_r,
$$
\textcolor{black}{with $\sK^{(s,x)}_r,\ \sB^{(s,x)}_r $ introduced 
\green in \eqref{SM1} and \eqref{SM2},}
and
\begin{align}
h(s,x,t,y):=\big(\nabla\tilde p^{(\tau,\xi)}(s,\cdot,t,y)(x)\big)_{(\tau,\xi)=(s,x)}
\stackrel{\eqref{KP8}}{=}-\nabla_{y} g^{(s,x)}_{s,t}(\theta_{s,t}(x)-y),\label{HF1}
\end{align}
and for a function $f$, 
$$
H_{s,t}f(x):=\int_{\mR^d}h(s,x,t,y)f(y)\dif y.
$$
By Lemma \ref{Le310} we have
$$
\nabla P_{s,t}f(x)=\nabla \wt P^{(\tau,\xi)}_{s,t}f(x)+\int^t_s\nabla \wt P^{(\tau,\xi)}_{s,r}\sA^{(\tau,\xi)}_r P_{r,t} f(x)\dif r.
$$
Taking $(\tau,\xi)=(s,x)$ and using the above notations, we can write
\begin{align}\label{AA87}
\nabla P_{s,t}f(x)=H_{s,t}f(x)+\int^t_sH_{s,r}\wt\sA_r P_{r,t} f(x)\dif r=H_{s,t}f(x)+\sum_{i=1}^4I^{(i)}_{s,t} f(x),
\end{align}
where for $u:=\frac{s+t}{2}$,
\begin{align*}
I^{(1)}_{s,t} f(x)&:=\int^{u} _sH_{s,r}\wt\sK_r P_{r,t}f(x)\dif r,\\
I^{(2)}_{s,t} f(x)&:=\int^{u} _sH_{s,r}\wt\sB_r P_{r,t}f(x)\dif r,\\
I^{(3)}_{s,t} f(x)&:=\int^t_{u}H_{s,r}\wt\sK_r P_{r,t}f(x)\dif r,\\
I^{(4)}_{s,t} f(x)&:=\int^t_{u}H_{s,r}\wt\sB_r P_{r,t}f(x)\dif r.
\end{align*}

{\it (Step 2).} Note that for $j\in\mN$, 
\begin{align}
|\nabla^j_yh(s,x,t,y)|&\stackrel{\eqref{HF1}}{=}|\nabla_{y}^{j+1} g^{(s,x)}_{s,t}(\theta_{s,t}(x)-y)|
\stackrel{\eqref{Es2}}{\lesssim}\phi^{(\alpha+j+1)}_{0,\alpha}(s,x,t,y).\label{DQ8}
\end{align}
Thus we have
\begin{align*}
|H_{s,t}f(x)|&\lesssim\int_{\mR^d}\phi^{(\alpha+1)}_{0,\alpha}(s,x,t,y)|f(y)|\dif y\\
&\leq (t-s)^{-\frac1\alpha}\int_{\mR^d}\phi^{(\alpha)}_{0,\alpha}(s,x,t,y)|f(y)|\dif y.
\end{align*}
For $I^{(1)}_{s,t} f(x)$, noting that by Lemma \ref{Le35},
\begin{align*}
|\wt\sK_r P_{r,t}f(z)|&\lesssim (1\wedge|z-\theta_{s,r}(x)|^{\gamma})|\cD^{(\alpha)}P_{r,t}f|(z)\\
&\lesssim (1\wedge|z-\theta_{s,r}(x)|^{\gamma})\int_{\mR^d}\phi^{(\alpha)}_{0,0}(r,z,t,y)|f(y)|\dif y,
\end{align*}
and using \eqref{AA1} and Lemma \ref{Le09}, we have
\begin{align*}
|I^{(1)}_{s,t} f(x)|&\stackrel{\eqref{DQ8}}{\lesssim}\int^{u}_s\!\!\!\int_{\mR^d}\phi^{(\alpha+1)}_{\gamma,\alpha}(s,x,r,z)
\int_{\mR^d}\phi^{(\alpha)}_{0,0}(r,z,t,y)|f(y)|\dif y\dif z\dif r\\
&\stackrel{\eqref{AA1}}{\lesssim}\int^{u}_s\int_{\mR^d}
\big(\phi^{(\alpha)}_{0,\alpha+\gamma-1}\odot\phi^{(\alpha)}_{0,0}\big)_r(s,x,t,y)|f(y)|\dif y\dif r\\
&\stackrel{\eqref{eq30}}{\lesssim}\int_{\mR^d}\phi^{(\alpha)}_{0,\alpha+\gamma-1}(s,x,t,y)|f(y)|\dif y.
\end{align*}
For $I^{(2)}_{s,t} f(x)$, noting that
$$
|\wt\sB_r P_{r,t}f(z)|\lesssim (|\theta_{s,r}(x)-z|^\beta+|\theta_{s,r}(x)-z|+(r-s)^{\beta/\alpha})|\nabla P_{r,t}f(z)|,
$$
using \eqref{DQ8} and \eqref{AA1}, we have
\begin{align*}
|I^{(2)}_{s,t} f(x)|&\lesssim\int^{u}_s\int_{\mR^d}\phi^{(\alpha)}_{0,\alpha+\beta-1}(s,x,r,z)|\nabla P_{r,t}f(z)|\dif z\dif r\\
&\lesssim(t-s)^{-\frac1\alpha}\int^t_s\int_{\mR^d}\phi^{(\alpha)}_{0,\alpha+\beta-1}(s,x,r,z)(t-r)^{\frac1\alpha}|\nabla P_{r,t}f(z)|\dif z\dif r.
\end{align*}

{\it (Step 3).} In this step we treat the hard term $I^{(3)}_{s,t} f(x)$. Let $\eps:=(t-r)^{1/\alpha}$ and
$$
\kappa_\eps(r,z,z'):=\kappa(r,\cdot,z')*\rho_\eps(z),\ \ 
\bar\kappa_\eps(r,z,z'):=\kappa_\eps(r,z,z')-\kappa_\eps(r,\theta_{r,s}(x),z')
$$
and
$$
\wt\sK^{(\eps)}_rf(z)=2\int_{\mR^d}\delta^{(1)}_f(z;z')\frac{\bar\kappa_\eps(r,z,z')}{|z'|^{d+\alpha}}\dif z'.
$$
Let us write
\begin{align*}
I^{(3)}_{s,t} f(x)&=\int^t_{u}\Big(H_{s,r}(\wt\sK_r-\wt\sK^{(\eps)}_r)P_{r,t}f(x)
+H_{s,r}\wt\sK^{(\eps)}_r P_{r,t}f(x)\Big)\dif r\\
&=:\int^t_{u}\Big(J^{(\eps)}_{1,r}(s,x,t)+J^{(\eps)}_{2,r}(s,x,t)\Big)\dif r.
\end{align*}
Let $\gamma_1\in(0,\gamma)$. Noting that
$$
|(\kappa-\kappa_\eps)(r,z,z')-(\kappa-\kappa_\eps)(r,\theta_{r,s}(x),z')|\lesssim_C(|z-\theta_{r,s}(x)|^{\gamma_1}\wedge 1)\eps^{\gamma-\gamma_1},
$$
by definition and Lemma \ref{Le35}, we have
\begin{align*}
|(\wt\sK_r-\wt\sK^{(\eps)}_r)P_{r,t}f(z)|&\lesssim 
(|z-\theta_{r,s}(x)|^{\gamma_1}\wedge 1)\eps^{\gamma-\gamma_1}|\cD^{(\alpha)}P_{r,t}f|(z)\\
&\lesssim (|z-\theta_{r,s}(x)|^{\gamma_1}\wedge 1)\eps^{\gamma-\gamma_1}\int_{\mR^d}\phi^{(\alpha)}_{0,0}(r,z,t,y)|f(y)|\dif y.
\end{align*}
For $J^{(\eps)}_{1,r}$, recalling $\eps=(t-r)^{1/\alpha}$,  we have
\begin{align*}
\int_{u}^t |J^{(\eps)}_{1,r}(s,x,t)|\dif r
&\stackrel{\eqref{DQ8}}{\lesssim}\int^t_{u}\int_{\mR^d}\big(\phi^{(\alpha)}_{0,\alpha+\gamma_1-1}\odot
\phi^{(\alpha)}_{0,\gamma-\gamma_1}\big)_r(s,x,t,y)|f(y)|\dif y\dif r\\
&\stackrel{\eqref{eq30}}{\lesssim}\int_{\mR^d}\phi^{(\alpha)}_{0,\alpha+\gamma-1}(s,x,t,y)|f(y)|\dif y.
\end{align*}
For $J^{(\eps)}_{2,r}$, by the change of variables and Fubini's theorem, we have
\begin{align*}
J^{(\eps)}_{2,r}(s,x,t)&=\int_{\mR^d}h(s,x,r,z)\int_{\mR^d}\delta^{(1)}_{P_{r,t}f}(z;z')\frac{\bar\kappa_\eps(r,z,z')}{|z'|^{d+\alpha}}\dif z'\dif z\\
&=\int_{\mR^d}\int_{\mR^d}\delta^{(1)}_{h(s,x,r,\cdot)\bar\kappa_\eps(r,\cdot,z')}(z;z')\frac{\dif z'}{|z'|^{d+\alpha}}P_{r,t}f(z)\dif z\\
&=\int_{\mR^d}h(s,x,r,z)\int_{\mR^d}\delta^{(1)}_{\bar\kappa_\eps(r,\cdot,z')}(z;z')\frac{\dif z'}{|z'|^{d+\alpha}}P_{r,t}f(z)\dif z\\
&+\int_{\mR^d}\int_{\mR^d}\delta^{(1)}_{h(s,x,r,\cdot)}(z;z')\bar\kappa_\eps(r,z+z',z')\frac{\dif z'}{|z'|^{d+\alpha}}P_{r,t}f(z)\dif z.
\end{align*}
Noting that by {\bf (H$^{\gamma}_a$)},
$$
|\delta^{(1)}_{\bar\kappa_\eps(r,\cdot,z')}(z;z')|\lesssim (\eps^{\gamma-1}|z'|)\wedge|z'|^\gamma\wedge 1,
$$
we have
\begin{align*}
\int_{\mR^d}|\delta^{(1)}_{\bar\kappa_\eps(r,\cdot,z')}(z;z')|\frac{\dif z'}{|z'|^{d+\alpha}}&\lesssim
\int_{\mR^d}((\eps^{\gamma-1}|z'|)\wedge|z'|^\gamma\wedge 1)\frac{\dif z'}{|z'|^{d+\alpha}}\lesssim\eps^{(\gamma-\alpha)\wedge 0}.
\end{align*}
On the other hand, by \eqref{HF1} and \eqref{Es2},
\begin{align*}
|\delta^{(1)}_{h(s,x,r,\cdot)}(z;z')|&\lesssim \big(((r-s)^{-\frac1\alpha}|z'|)\wedge 1\big)
\big(\phi^{(\alpha)}_{0,\alpha-1}(s,x,r,z+z')+\phi^{(\alpha)}_{0,\alpha-1}(s,x,r,z)\big).
\end{align*}
Thus, as in \eqref{AA86} we have
\begin{align*}
\int_{\mR^d}\frac{|\delta^{(1)}_{h(s,x,r,\cdot)}(z;z')|\dif z'}{|z'|^{d+\alpha}}
&\lesssim \int_{\mR^d}\big(((r-s)^{-\frac1\alpha}|z'|)\wedge 1\big)\phi^{(\alpha)}_{0,\alpha-1}(s,x,r,z+z')\frac{\dif z'}{|z'|^{d+\alpha}}\\
&\quad+\phi^{(\alpha)}_{0,\alpha-1}(s,x,r,z)\int_{\mR^d}\big(((r-s)^{-\frac1\alpha}|z'|)\wedge 1\big)\frac{\dif z'}{|z'|^{d+\alpha}}\\
&\lesssim \phi^{(\alpha)}_{0,\alpha-1}(s,x,r,z)(r-s)^{-1}=\phi^{(\alpha)}_{0,-1}(s,x,r,z).
\end{align*}
Therefore,
\begin{align*}
|J^{(\eps)}_{2,r}(s,x,t)|
&\lesssim \int_{\mR^d}\Big[\eps^{
(\gamma-\alpha)\wedge 0
}\phi^{(\alpha)}_{0,\alpha-1}(s,x,r,z)
+\phi^{(\alpha)}_{0,-1}(s,x,r,z)\Big]P_{r,t}|f|(z)\dif z.
\end{align*}
Recall $\eps=(t-r)^{\frac1\alpha}$. By \eqref{AA1}, we obtain
$$
\int_{u}^t |J^{(\eps)}_{2,r}(s,x,t)|\dif r\lesssim(t-s)^{-\frac1\alpha}\int_{\mR^d}\phi^{(\alpha)}_{0,\alpha}(s,x,t,y)|f(y)|\dif y.
$$

{\it (Step 4).} For $\eps=(t-r)^{1/\alpha}$, we define
$$
\bar b_\eps(r,z):=(b*\rho_\eps)(r,z)-(b*\rho_{\eps}*\rho_{|r-s|^{1/\alpha}})(r,\theta_{r,s}(x))
$$
and
$$
\wt\sB^{(\eps)}_r f(z):=\bar b_\eps(r,z)\cdot\nabla f(z).
$$
For $I^{(4)}_{s,t} $, we similarly write
\begin{align*}
I^{(4)}_{s,t} f(x)&=\int^t_{u}\Big(H_{s,r}(\wt\sB_r-\wt\sB^{(\eps)}_r)P_{r,t}f(x)
+H_{s,r}\wt\sB^{(\eps)}_r P_{r,t}f(x)\Big)\dif r\\
&=:\int^t_{u}\Big(J^{(\eps)}_{3,r}(s,x,t)+J^{(\eps)}_{4,r}(s,x,t)\Big)\dif r.
\end{align*}
For $J^{(\eps)}_{3,r}$, since
$$
|\bar b_0-\bar b_\eps|(r,z)\leq \kappa_0\eps^\beta=\kappa_0(t-s)^{\beta/\alpha},
$$
by \eqref{DQ8} we have
\begin{align*}
|J^{(\eps)}_{3,r}(s,x,t)|&=\left|\int_{\mR^d}h(s,x,r,z)(\bar b_0(r,z)-\bar b_\eps(r,z))\cdot \nabla P_{r,t}f(z)\dif z\right|\\
&\lesssim\int_{\mR^d}\phi^{(\alpha+1)}_{0,\alpha}(s,x,r,z)(t-r)^{\frac{\beta}{\alpha}}|\nabla P_{r,t}f(z)|\dif z,
\end{align*}
and
$$
\int^t_{u}|J^{(\eps)}_{\textcolor{black}{3},r}(s,x,t)|\dif r\lesssim 
(t-s)^{-\frac{1}{\alpha}}\int^t_s\!\!\int_{\mR^d}\phi^{(\alpha)}_{0,\alpha}(s,x,r,z)(t-r)^{\frac{\beta}{\alpha}}|\nabla P_{r,t}f(z)|\dif z\dif r.
$$
For $J^{(\eps)}_{4,r}$, we derive integrating by parts that
\begin{align*}
|J^{(\eps)}_{4,r}(s,x,t)|&=\left|\int_{\mR^d}h(s,x,r,z)\, \bar b_\eps(r,z)\cdot\nabla_z P_{r,t}f(z)\dif z\right|
\\&\leq\left|\int_{\mR^d}h(s,x,r,z)\, \div \bar b_\eps(r,z) P_{r,t}f(z)\dif z\right|
\\&\quad+\left|\int_{\mR^d} \bar b_\eps(r,z)\cdot \nabla_z h(s,x,r,z) P_{r,t}f(z)\dif z\right|.
\end{align*}
Since
$$
|\div \bar b_\eps(r,z)|=|\div  b_\eps(r,z)|\leq \kappa_0\eps^{\beta-1}=\kappa_0(t-r)^{(\beta-1)/\alpha}
$$
and
$$
|\bar b_\eps|(r,z)\lesssim |z-\theta_{r,s}(x)|^\beta+(r-s)^{\beta/\alpha},
$$
by \eqref{DQ8} and \eqref{AA1} we have
\begin{align*}
|J^{(\eps)}_{4,r}(s,x,t)|&\lesssim\int_{\mR^d}\phi^{(\alpha)}_{0,\alpha-1}(s,x,r,z) (t-r)^{(\beta-1)/\alpha}|P_{r,t}f(z)|\dif z\\
&\quad+\int_{\mR^d}\phi^{(\alpha)}_{0,\alpha+\beta-2}(s,x,r,z)|P_{r,t}f(z)|\dif z.
\end{align*}
Thus,  
\begin{align*}
\int^t_{u}|J^{(\eps)}_{4,r}(s,x,t)|\dif r&\lesssim\int^t_{u}\int_{\mR^d}\big(\phi^{(\alpha)}_{0,\alpha-1}
\odot\phi^{(\alpha)}_{0,\alpha+\beta-1}\big)_r(s,x,t,y)|f(y)|\dif y\dif r\\
&\qquad+\int^t_{u}\int_{\mR^d}\big(\phi^{(\alpha)}_{0,\alpha+\beta-2}
\odot\phi^{(\alpha)}_{0,\alpha}\big)_r(s,x,t,y)|f(y)|\dif y\dif r\\
&\stackrel{\eqref{eq30}}{\lesssim}\int^t_u\big[(r-s)^{\frac{\alpha-1}{\alpha}}(t-r)^{\frac{\alpha+\beta-1}{\alpha}}+(r-s)^{\frac{\alpha+\beta-2}{\alpha}}(t-r)\big]\\
&\qquad\times\big[(r-s)^{-1}+(t-r)^{-1}\big]\dif r
\int_{\mR^d}\phi^{(\alpha)}_{0,0}(s,x,t,y)|f(y)|\dif y\\
&\lesssim\int_{\mR^d}\phi^{(\alpha)}_{0,2\alpha+\beta-2}(s,x,t,y)|f(y)|\dif y\\
&\lesssim (t-s)^{-\frac 1\alpha}\int_{\mR^d}\phi^{(\alpha)}_{0,\alpha}(s,x,t,y)|f(y)|\dif y,
\end{align*}
recalling that $\alpha+\beta>1 $ for the last inequality.
Hence,
\begin{align*}
|I^{(4)}_{s,t} f(x)|
&\lesssim(t-s)^{-\frac{1}{\alpha}}\Big( \int^t_s\int_{\mR^d}\phi^{(\alpha)}_{0,\alpha}(s,x,r,z)(t-r)^{\frac{\beta}{\alpha}}|\nabla P_{r,t}f(z)|\dif z\dif r\\
&\quad+\int_{\mR^d}\phi^{(\alpha)}_{0,\alpha}(s,x,t,y)|f(y)|\dif y\Big).
\end{align*}

{\it (Step 5).} Combining the above calculations, we obtain \label{USE_THM_4.3}
\begin{align*}
&|\nabla P_{s,t} f(x)|\lesssim(t-s)^{-\frac1\alpha}\int_{\mR^d}\phi^{(\alpha)}_{0,\alpha}(s,x,t,y)|f(y)|\dif y\\
&+(t-s)^{-\frac{1}{\alpha}}\int^t_s\int_{\mR^d}\phi^{(\alpha)}_{0,\alpha}(s,x,r,z)(t-r)^{\frac{\beta}{\alpha}}|\nabla P_{r,t}f(z)|\dif z\dif r\\
&+(t-s)^{-\frac{1}{\alpha}}\int^t_s\int_{\mR^d}\phi^{(\alpha)}_{0,\alpha+\beta-1}(s,x,r,z)(t-r)^{\frac1\alpha}|\nabla P_{r,t}f(z)|\dif z\dif r.
\end{align*}
By the lower bound estimate, we further have
\begin{align}\label{DQ2}
\begin{split}
(t-s)^{\frac1\alpha}|\nabla P_{s,t} f(x)|&\lesssim P_{s,t}|f|(x)+\int^t_s(t-r)^{\frac\beta\alpha}P_{s,r}|\nabla P_{r,t}f|(x)\dif r\\
&\quad+\int^t_s(r-s)^{\frac{\beta-1}{\alpha}}(t-r)^{\frac1\alpha}P_{s,r}|\nabla P_{r,t}f|(x)\dif r.
\end{split}
\end{align}
For fixed $0\leq u<t\leq T$ and $s\in(u,t)$, we let
$$
\Gamma^t_u(s,x):=(t-s)^{\frac1\alpha}P_{u,s}|\nabla P_{s,t}f|(x).
$$
Using $P_{u,s}$ act on both sides of \eqref{DQ2} and by $P_{u,s}P_{s,r}=P_{u,r}$, we derive that
\begin{align*}
\Gamma^t_u(s,x)&\lesssim P_{u,t}|f|(x)
+\int^t_s\Big[(r-s)^{\frac{\beta-1}{\alpha}}+(t-r)^{\frac{\beta-1}{\alpha}}\Big]\Gamma^t_u(r,x)\dif r.
\end{align*}
Note that by Theorem \ref{Le36},
$$
\sup_{s\in[u,t]}\|\Gamma^t_u(s,\cdot)\|_\infty<\infty.
$$
Since $\alpha+\beta>1$, from the Volterra-Gronwall inequality, we obtain that for all $s\in(u,t)$,
$$
\Gamma^t_u(s,x)\lesssim P_{u,t}|f|(x).
$$
Taking limit $u\uparrow s$, we obtain
$$
(t-s)^{\frac1\alpha}|\nabla P_{s,t}f|(x)\lesssim P_{s,t}|f|(x),
$$
which eventually yields the desired gradient estimate.

{\it (Step 6).}  Finally, by \eqref{AA87} and the dominated convergence theorem, one can show that
$$
\lim_{x\to x_0}\sup_{(b,a)\in\sC}|\nabla P^{b,a}_{s,t}f(x)-\nabla P^{b,a}_{s,t}f(x_0)|=0.
$$
Indeed, from the above proof, it suffices to show that
$$
\lim_{x\to x_0}\sup_{(b,a)\in\sC}|H^{b,a}_{s,t}f(x)-H^{b,a}_{s,t}f(x_0)|=0.
$$
This follows by Lemma \ref{Le27}.

\section{Proof of Theorem \ref{Main}}
\label{FINAL_SEC}

The point is here to prove Theorem \ref{Main}. Namely, we want to extend the bounds of Theorem \ref{T8}, Lemma \ref{Le35} and Theorem \ref{Th41} under the sole assumptions {\bf (H$^\gamma_a$)}, {\bf (H$^\beta_b$)}.
{\blue In Subsection \ref{SUB_SEC_BD} we consider bounded drift $b$ (without additional smoothness than in {\bf (H$^\gamma_a$)}, {\bf (H$^\beta_b$)}). 
In Subsection \ref{SUB_SEC_TRUNC} we drop the bounded assumption on $b$ by a truncation argument.}
\subsection{Bounded drift $b$}\label{SUB_SEC_BD}
{\blue In this subsection we assume  {\bf (H$^\gamma_a$)}, {\bf (H$^\beta_b$)} and $b$ is bounded.}
Let $a_\eps$ and $b_\eps$ be the smooth approximations of $a$ and $b$, respectively. Hence, assumptions {\bf (H$^\gamma_a$)}, {\bf (H$^\beta_b$)} and \eqref{DA1} are met by $a_\eps, b_\eps $ for the SDE
\begin{equation}
\dif X_t^\eps=b_\eps(t,X_t^\eps)\dif t+a_\eps(t,X_{t-}^\eps)\dif L^{(\alpha)}_t.
\label{SDE_WITH_MOLL}
\end{equation} 
The following convergence in law result was established in \cite{CHZ}, see  Theorem 1.1 therein.
\bt\label{Th62}
Let $X^\eps_{s,t}(x)$ be the unique solution of SDE \eqref{SDE_WITH_MOLL} starting from $x$ at time $s$. Then $X^\eps_{s,t}(x)$ weakly converges to $X_{s,t}(x)$.
\et
\begin{proof}
For fixed $(s,x)\in\mR_+\times\mR^d$, since the coefficients $b,a$ \textcolor{black}{have} linear growth, it is by now standard to show that the law of 
$X^\eps_{s,\cdot}(x)$ is tight in the space of all c\`adl\`ag  functions (\textcolor{black}{see e.g. the Aldous criterion in the monographs \cite[Chaper VI, Theorem 4.5]{jaco:shir:87}, \cite{bass:11}}).  \textcolor{black}{Through the martingale problem formulation}, one can show \textcolor{black}{as well} that any weak accumulation point of the law of
$X^\eps_{s,\cdot}(x)$ is a weak solution of SDE \eqref{SDE0}. Finally, by the weak uniqueness, one sees that 
$X^\eps_{s,t}(x)$ weakly converges to $X_{s,t}(x)$.
\end{proof}

Denoting by $p_\eps $ the associated density, it therefore holds from Theorem \ref{T8},  
Lemma \ref{Le35} and Theorem \ref{Th41} that
\begin{enumerate}[(i)]
\item {\bf (Two-sides estimate)} For any $T>0$, there is a constant $C_1=C_1(T,\Theta)>0$ such that for all $0\leq s<t\leq T$ and $x,y\in\mR^d$,
\begin{align}\label{TW_EPS}
p_\eps(s,x,t,y)\asymp_{C_1} \phi^{(\alpha)}_{0,\alpha}(s,x,t,y).
\end{align}
\item {\bf (Fractional derivative estimate)} For any $T>0$, there is a constant $C_2=C_2(T,\Theta)>0$ such that for all $0\leq s<t\leq T$ and $x,y\in\mR^d$,
\begin{align}\label{FR_EPS}
|\cD^{(\alpha)}{p_\eps(s,\cdot,t,y)}|(x)\lesssim_{C_2}\phi^{(\alpha)}_{0,0}(s,x,t,y).
\end{align}
\item {\bf (Gradient estimate in $x$)} For any $T>0$, there is a constant $C_3=C_3(T,\Theta)>0$ such that for all $0\leq s<t\leq T$ and $x,y\in\mR^d$,
\begin{align}\label{GR_EPS}
|\nabla P^\eps_{s,t}f(x)|\lesssim_{C_3} (t-s)^{-1/\alpha}P^\eps_{s,t}|f|(x).
\end{align}
\end{enumerate}
where the constants in the above controls only depend on {\bf (H$^\gamma_a$)}, {\bf (H$^\beta_b$)} through $\Theta $ 
(see precisely \eqref{DEF_THETA}).

By Theorem \ref{Th62}, we have for any $f\in C_b(\mR^d)$,
\begin{align}\label{LL9}
\lim_{\eps\to 0}P^\eps_{s,t}f(x):=\lim_{\eps\to 0}\mE f(X^\eps_{s,t}(x))=\mE f(X_{s,t}(x))=:P_{s,t}f(x).
\end{align}

(i) (Two-sided estimates) For nonnegative measurable function $f$, we get from \eqref{TW_EPS} 
\begin{align*}
C^{-1}_1\int_{\R^d}\phi^{(\alpha)}_{0,\alpha}(s,x,t,y)f(y)\dif y\leq \mE f(X_{s,t}(x))\leq C_1\int_{\R^d}\phi^{(\alpha)}_{0,\alpha}(s,x,t,y) f(y)\dif y,
\end{align*}
which implies that $X_{s,t}(x)$ has a density $p(s,x,t,y)$ having lower and upper bound as in \eqref{TW}. 
On the other hand, for fixed $s<t$, by Theorem \ref{T8} we have
$$
(x,y)\mapsto p_\eps(s,x,t,y) \mbox{ is equi-continuous in $\eps\in(0,1)$}.
$$
From the Ascoli-Arzel\`a theorem, there are a subsequence $\eps_k$ and a continuous function $\bar p(s,x,t,y)$ as a function of $x,y\in\mR^d$ such that 
\begin{align}\label{ED8}
p_{\eps_k}(s,x,t,y)\to \bar p(s,x,t,y) \mbox{ locally uniformly in $x,y\in\mR^d$},
\end{align}
which together with \eqref{LL9} yields that
\begin{align}\label{ED9}
p(s,x,t,y)=\bar p(s,x,t,y)\mbox{ is continuous as a function of $x,y\in\mR^d$.}
\end{align}

(ii) (Fractional derivative estimates)
It follows by \eqref{FR_EPS}, \eqref{ED8}, \eqref{ED9} and Fatou's lemma that
\begin{align*}
|\cD^{(\alpha)}{p(s,\cdot,t,y)}|(x)&=\int_{\mR^d}\lim_{k\to\infty}|\delta^{(2)}_{p_{\eps_k}(s,\cdot,t,y)}(x;z)|\frac{\dif z}{|z|^{d+\alpha}}\\
&\leq\varliminf_{k\to\infty}\int_{\mR^d}|\delta^{(2)}_{p_{\eps_k}(s,\cdot,t,y)}(x;z)|\frac{\dif z}{|z|^{d+\alpha}}\\
&=\varliminf_{k\to\infty}|\cD^{(\alpha)}{p_{\eps_k}(s,\cdot,t,y)}|(x)\lesssim_{C_2}\phi^{(\alpha)}_{0,0}(s,x,t,y).
\end{align*}

(iii) (Gradient estimates) For fixed $f\in C_b(\mR^d)$, by \eqref{GR_EPS}, 
$$
x\mapsto\nabla P^\eps_{s,t}f(x) \text{ is equi-continuous in $\eps$},
$$
which together with \eqref{LL9} implies that $x\mapsto P_{s,t}f(x)$ is continuous differentiable. 
By taking limits along a subsequence $\eps_k$ for \eqref{GR_EPS}, we obtain
$$
|\nabla P_{s,t}f(x)|\lesssim_{C_3} (t-s)^{-1/\alpha}P_{s,t}|f|(x).
$$
Finally, for fixed $t'>t$ and $y\in\mR^d$, we let $f(x):=p(t,x,t',y)$, then by the Chapman-Kolmogorov equation, we obtain
$$
|\nabla p(s,\cdot,t',y)(x)|\lesssim_{C_3} (t-s)^{-1/\alpha}p(s,x,t',y).
$$
This then readily gives estimate \eqref{GR} (logarithmic derivative) of Theorem \ref{Main}.\\

\subsection{Unbounded drift $b$} \label{SUB_SEC_TRUNC}
{\blue In this subsection we assume  {\bf (H$^\gamma_a$)} and {\bf (H$^\beta_b$)}.
For $n\in\mN$, define
$$
b^n(t,x):=(-n)\vee b(t,x)\wedge n.
$$
With the $\kappa_0$ as in \eqref{BB0}, we have
 \begin{align}\label{BB00}
|b^n(t,0)|\le \kappa_0,\ |b^n(t,x)-b^n(t,y)|\leq\kappa_0(|x-y|^\beta\vee|x-y|).
\end{align}
Consider the following SDE
\begin{equation}
\dif X_t^n=b^n(t,X_t^n)\dif t+a(t,X_{t-}^n)\dif L^{(\alpha)}_t.
\label{SDE_WITH_TRUN}
\end{equation} 
The following result is the very similar to Theorem \ref{Th62}.
\bt\label{Th622}
Let $X^n_{s,t}(x)$ be the unique solution of SDE \eqref{SDE_WITH_TRUN} starting from $x$ at time $s$. Then $X^n_{s,t}(x)$ weakly converges to $X_{s,t}(x)$.
\et
Moreover, for fixed $s>0$, consider the following ODE:
\begin{align}\label{ODE10}
\dot\theta^n_{s,t}=b^n_{|t-s|^{1/\alpha}}(t,\theta^n_{s,t}),\ \ \theta^n_{s,s}=x,\ \ \ t\geq 0.
\end{align}
We have the following convergence.
\bl\label{Le53}
For each $s,t>0$ and $x\in\mR^d$, it holds that
$$
\lim_{n\to\infty}|\theta^n_{s,t}(x)-\theta_{s,t}(x)|=0.
$$
\el
\begin{proof}
We drop the starting point $x$ and assume $t>s$. From dynamics \eqref{ODE1} and \eqref{ODE10},
\begin{align*}
|\theta^n_{s,t}-\theta_{s,t}|
&\leq\int^t_s\|\nabla b^n_{|r-s|^{1/\alpha}}(r,\cdot)\|_\infty|\theta^n_{s,r}-\theta_{s,r}|\dif r\\
&+\int^t_s\big|b^n_{|r-s|^{1/\alpha}}(r,\theta_{s,r})-b_{|r-s|^{1/\alpha}}(r,\theta_{s,r})\big|\dif r.
\end{align*}
By \eqref{BB00}, \eqref{LL01} and Gronwall's inequality, we have
$$
|\theta^n_{s,t}-\theta_{s,t}|\leq C\int^t_s\big|b^n_{|r-s|^{1/\alpha}}(r,\theta_{s,r})-b_{|r-s|^{1/\alpha}}(r,\theta_{s,r})\big|\dif r,
$$
which gives the desired limit by the dominated convergence theorem.
\end{proof}
Now by Theorem \ref{Th622} and Lemma \ref{Le53}, using exactly the  same argument as in Subsection \ref{SUB_SEC_BD}, we can show Theorem \ref{Main}.
}

\medskip

{\blue {\bf Acknowledgment:} 
The authors would like to thank the referee for his/her very useful comments.
}

\end{document}